\documentclass[11pt,reqno]{amsart}
\usepackage{latexsym,amssymb,amsfonts,amsmath,graphicx,url}
\usepackage[margin=1in]{geometry}
\usepackage{stmaryrd}
\usepackage[vcentermath,enableskew]{youngtab}
\usepackage{color}
\usepackage[all]{xy}
\usepackage{makecell}
\usepackage{gensymb}
\usepackage{caption}
\usepackage{subcaption}
\usepackage{ytableau}
\usepackage{cases}
\usepackage{mathtools}

\newtheorem{theorem}{Theorem}[section]
\newtheorem{proposition}[theorem]{Proposition}
\newtheorem{lemma}[theorem]{Lemma}
\newtheorem{conjecture}[theorem]{Conjecture}
\newtheorem{corollary}[theorem]{Corollary}

\theoremstyle{definition}
\newtheorem{definition}[theorem]{Definition}
\newtheorem{example}[theorem]{Example}

\theoremstyle{remark}
\newtheorem{remark}[theorem]{Remark}

\newcommand{\orb}{\mathcal{O}}

\DeclareMathOperator*{\pro}{Pro}
\DeclareMathOperator*{\row}{Row}

\DeclareMathOperator*{\togpro}{TogPro}

\usepackage{tikz}
\DeclareRobustCommand{\widetriangle}{%
\begin{tikzpicture}%
\draw (-1.5ex,0) -- (1.5ex, 0) -- (0, 2ex) -- (-1.5ex,0);
\end{tikzpicture}%
}

\DeclareRobustCommand{\subwidetriangle}{%
\begin{tikzpicture}%
\draw (-.75ex,0) -- (.75ex, 0) -- (0, 1ex) -- (-.75ex,0);
\end{tikzpicture}%
}

\title[$P$-strict promotion and $Q$-partition rowmotion]{$P$-strict promotion and $Q$-partition rowmotion: \linebreak the graded case}
\author[Bernstein, Striker, Vorland]{Joseph Bernstein, Jessica Striker, and Corey Vorland}
\email{joseph.bernstein@ndsu.edu, jessica.striker@ndsu.edu, corey.vorland@gmail.com}

\address{North Dakota State University}

\begin{document}
\maketitle

\begin{abstract}
  Promotion and rowmotion are intriguing actions in dynamical algebraic combinatorics which have inspired much work in recent years. In this paper, we study $P$-strict labelings of a finite, graded poset $P$ of rank $n$ and labels at most $q$, which generalize semistandard Young tableaux with $n$ rows and entries at most $q$, under promotion. These $P$-strict labelings are in equivariant bijection with $Q$-partitions under rowmotion, where $Q$ equals the product of $P$ and a chain of $q-n-1$ elements. We study the case where $P$ equals the product of chains in detail, yielding new homomesy and order results in the realm of tableaux and beyond. Furthermore, we apply the bijection to the cases in which $P$ is a minuscule poset and when $P$ is the three element $V$ poset. Finally, we give resonance results for promotion on $P$-strict labelings and rowmotion on $Q$-partitions.
\end{abstract}

\tableofcontents

\section{Introduction}
\emph{Semistandard Young tableaux} and \emph{Gelfand--Tsetlin patterns} are well-loved combinatorial objects with a nice, statistic-preserving bijection between them. In \cite{BSV2021}, we generalized this correspondence to objects we called \emph{$P$-strict labelings of} $P\times[\ell]$ (analogous to semistandard tableaux) for a finite poset $P$ and $\ell\in\mathbb{N}$, and \emph{$B$-bounded $Q$-partitions} (analogous to Gelfand--Tsetlin patterns) for related posets $Q$ and bounding function $B$. In addition, we showed this bijection is \emph{equivariant}, mapping the well-studied \emph{promotion} action on tableaux to a \emph{piecewise-linear toggle group action} on $Q$-partitions, which under a certain condition is equivalent to \emph{rowmotion}, an important action in dynamical algebraic combinatorics.

Informally, $P$-strict labelings of $P\times[\ell]$ are labelings of $P\times[\ell]$ with positive integers that strictly increase on each copy of $P$ and weakly increase along each copy of $[\ell]$. Additional parameters include a \emph{restriction function} $R$, that specifies which labels are allowed in the $P$-strict labeling, as well as functional parameters $u$ and $v$. The case $P=[n]$ corresponds to (skew) semistandard Young tableaux of $n$ rows, where $\ell$ is the number of tableau columns, and $u$ and $v$ determine the shape of the skew tableau by specifying which partitions to remove from the upper left and lower right of the $n\times\ell$ bounding rectangle (the case of $u=v=0$ corresponds to rectangular tableaux). While the bijection of \cite{BSV2021} allows for all of these parameters in full generality, our applications of the general formula in that paper were in the case $P=[n]$. These included interesting results and conjectures on flagged and symplectic tableaux.

In this paper, we apply the results of \cite{BSV2021} to cases of interest beyond $P=[n]$, with an eye toward translating objects with known dynamical behavior through our bijection in order to obtain new results and new perspectives on conjectures. Along with periodicity, we aim to describe the behavior of promotion and rowmotion through the notions of homomesy and resonance.  The \emph{homomesy phenomenon}, first described in \cite{PR2015}, is where the average value of a given combinatorial statistic is the same across all orbits.  For an overview of this phenomenon, including many examples, see the survey of T.~Roby \cite{Roby2016}.  \emph{Resonance}, defined in \cite{DPS2017}, is a way to capture the structure of orbits without requiring predictable periodicity, by projecting a more complicated object and action to a simpler set that demonstrates cyclic behavior.  The article \cite{Striker2017} contains a more detailed description of resonance, as well as a general survey of dynamical algebraic combinatorics, including homomesy.

This paper is organized as follows. In Section~\ref{sec:gen}, we give relevant background definitions and theorems from \cite{BSV2021} at the level of generality necessary for this paper ($u=v=0$) and state Corollaries~\ref{thm:PrankedGlobalq} and \ref{cor:PrankedGlobalqRow}, which are consequences of these general results when $P$ is graded. In Section~\ref{sec:ab}, we apply these results to the case where $P$ is a product of chains, proving distribution results on both sides of the bijection (Theorems~\ref{thm:ppartdist} and \ref{thm:pstrictdist}). Along the way, we also discover a new distribution result on semistandard Young tableaux under promotion (Theorem~\ref{ssytbcdist}). Section~\ref{sec:beyond} considers the dynamics of $P$-strict labelings in cases where either $P$ is graded but not equal to a product of chains or the labels of $P \times [\ell]$ are restricted by \emph{flags} instead of by a global bound. Finally, in Section~\ref{sec:resonance}, we prove Theorem~\ref{thm:ConResonance}, a resonance result on general $P$-strict labelings with global bound $q$.

\section{Background and general results}
\label{sec:gen}
In this section, we give relevant background and general results from \cite{BSV2021}.
In Subsections~\ref{sec:pro}-\ref{sec:results}, we state definitions and theorems in the  less general case of $u=v=0$ and restriction function $R_{\alpha}^{\beta}$. Then in Subsection~\ref{sec:newcor}, we restrict our attention to the case where $P$ is graded and state  consequences of these theorems that will be of use in the remainder of the paper. For the definitions and theorems in full generality, see \cite{BSV2021}.

\subsection{Promotion on $P$-strict labelings}
\label{sec:pro}
The definitions in this section are adapted from \cite[\textsection 1.2]{BSV2021}. We begin with the following preliminary definitions.
\begin{definition} In this paper, {$P$ and $Q$ represent finite posets} with partial orders $\leq_P$ and $\leq_Q$, respectively. Also, $\lessdot$ indicates a {covering relation} in a poset, $J(P)$ is the poset of order ideals of $P$ ordered by containment, {$\ell$ and $q$ are positive integers}, $[\ell]$ denotes a {chain poset} of $\ell$ elements, and $P\times[\ell]=\{(p,i) \ | \ p\in P, i\in \mathbb{N}, 1\leq i\leq \ell\}$. $\mathcal{P}(\mathbb{Z})$ represents the set of all nonempty, finite subsets of $\mathbb{Z}$.

  Let $L_i$ represent the $i$th \textbf{layer} of $P\times[\ell]$, that is, the set of $(p,i)\in P\times[\ell]$ where $p$ ranges over all possible values and $i$ is fixed. Let $F_p$ denote the $p$th \textbf{fiber} of $P\times[\ell]$, that is, the set of $(p,i)\in P\times[\ell]$ where $i$ ranges over all possible values and $p$ is fixed.
\end{definition}

\begin{definition}
  \label{def:Pstrict}
  We say that a function $f:P \times [\ell] \rightarrow \mathbb{Z}^+$ is a $P$-\textbf{strict labeling of $P \times [\ell]$ with restriction function $R_{\alpha}^{\beta}:P \rightarrow \mathcal{P}(\mathbb{Z}^+)$}, $\alpha, \beta: P \to \mathbb{Z}^+$ if $f$ satisfies the following:
  \begin{enumerate}
    \item $f(p_1,i)<f(p_2,i)$ whenever $p_1<_P p_2$,
    \item $f(p,i_1)\leq f(p,i_2)$ whenever $i_1\leq i_2$,
    \item $\alpha(p) \leq f(p,i) \leq \beta(p)$.
  \end{enumerate}
  That is, $f$ is strictly increasing inside each copy of $P$ (layer), weakly increasing along each copy of the chain $[\ell]$ (fiber), and such that $\alpha(p)$ and $\beta(p)$ are the lower and upper bounds, respectively, of the labels on the $p$th fiber.
  Let $\mathcal{L}_{P \times [\ell]}(R_{\alpha}^{\beta})$ denote the set of all $P$-strict labelings on $P \times [\ell]$ with restriction function $R_{\alpha}^{\beta}$.
\end{definition}

 In the case of a global upper bound $q$, our restriction function will be $R^q_1$, that is, we take $\alpha$ to be the constant function $1$ and $\beta$ to be the constant function $q$. Since a lower bound of $1$ is used frequently, we suppress the subscript $1$; that is, if no subscript appears, we take it to be $1$.

\begin{example}
  The objects on the left half of Figure~\ref{fig:PrankedGlobalq} are elements of $\mathcal{L}_{P \times [4]}(R^6)$, where $P$ is the X-shaped poset shown in the center. The labels on each of the four layers are strictly increasing (from $1$ up to the global maximum $q = 6$) and are connected by solid lines, while the labels along each of the five fibers are weakly increasing and connected by dotted lines.
\end{example}

\begin{definition}
  We say a label $f(p,i)$ in a $P$-strict labeling $f\in \mathcal{L}_{P \times [\ell]}(R_\alpha^\beta)$ is \textbf{raisable (lowerable)} if there exists another $P$-strict labeling  $g\in\mathcal{L}_{P \times [\ell]}(R_\alpha^\beta)$ where $f(p,i)<g(p,i)$ ($f(p,i)>g(p,i)$), and $f(p',i')=g(p',i')$ for all $(p',i')\in P \times [\ell]$, $p' \neq p$.
\end{definition}

\begin{definition}
  \label{def:BenderKnuth}
  Let the action of the \textbf{$k$th Bender--Knuth involution} $\rho_k$ on
  a $P$-strict labeling $f\in \mathcal{L}_{P \times [\ell]}(R_\alpha^\beta)$ be as follows:
  identify all raisable labels $f(p,i)=k$ and all lowerable labels $f(p,i)=k+1$.
  Call these labels `free'. Suppose the labels $f(F_p)$ include $a$ free $k$ labels followed by $b$ free $k+1$ labels; $\rho_k$ changes these labels to $b$ copies of $k$ followed by $a$ copies of $k+1$.
  \textbf{Promotion} on $P$-strict labelings is defined as the composition of these involutions: $\pro(f)=\cdots\circ\rho_3\circ\rho_2\circ\rho_1(f)$. Note that since $\beta$ induces upper bounds on the labels, only a finite number of Bender--Knuth involutions act nontrivially.
\end{definition}

\begin{figure}[hbtp]
  \begin{center}
    \includegraphics[width=\textwidth]{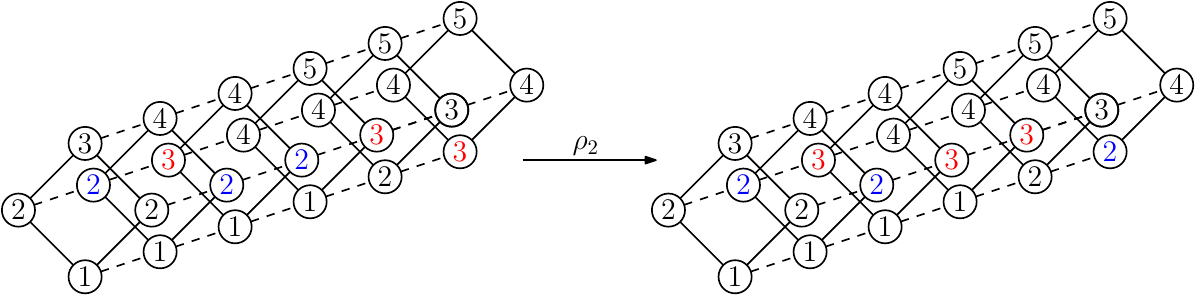}
  \end{center}
  \caption{An example of the Bender--Knuth involution $\rho_2$ on an element of $\mathcal{L}_{P \times [6]}(R^5)$, where $P$ is the four element diamond shaped poset. The raisable 2 labels are colored blue and the lowerable 3 labels are colored red.} \label{fig:DiamondBK}
\end{figure}

See Figure~\ref{fig:DiamondBK} for an example of a Bender--Knuth involution on a $P$-strict labeling.

\subsection{Rowmotion on $Q$-partitions}
\label{sec:row}
The definitions in this section are adapted from \cite[\textsection 1.3]{BSV2021}.
\begin{definition}
  \label{def:ppart}
  Let $Q$ be a poset. A \textbf{$Q$-partition} is a map $\sigma:Q\rightarrow\mathbb{Z}_{\geq 0}$ such that if $x\leq_Q x'$, then $\sigma(x)\leq \sigma(x')$. Let $\widehat{Q}$ denote $Q$ with $\widehat{0}$ added below all elements and $\widehat{1}$ added above all elements. Let $\mathcal{A}^{\ell}({Q})$ denote the set of all $\widehat{Q}$-partitions $\sigma$ with $\sigma(\widehat{0})=0$ and $\sigma(\widehat{1})=\ell$.
\end{definition}
Note we use $Q$ for a generic poset in this context rather than $P$ to avoid confusion when we later relate these objects to the objects of the previous subsection.

\begin{example}
  The objects on the right half of Figure \ref{fig:PrankedGlobalq} are elements of $\mathcal{A}^{4}(P \times [3])$, where $P$ is the X-shaped poset shown in the center.  In our visualizations, we omit the elements $\widehat{0}$ and $\widehat{1}$.
\end{example}

\begin{remark}
  \label{rem:JofP}
  When $\ell=1$, $\mathcal{A}^{\ell}({Q})=J(Q)$, the set of {order ideals (or lower sets)} of $Q$. The set of order ideals forms a \emph{distributive lattice} ordered by containment. This is the setting in which rowmotion and toggles were originally studied~\cite{CF1995,SW2012}.
\end{remark}

In Definitions~\ref{def:toggle1} and \ref{def:row1} below, we define toggles and rowmotion in what is often called the \emph{piecewise-linear} context. These definitions are equivalent (by rescaling) to those first given by Einstein and Propp on the order polytope~\cite{EP2014,EP2021}.

\begin{definition}
  \label{def:toggle1}
  For $\sigma \in \mathcal{A}^{\ell}({Q})$ and $x \in Q$, let $\nabla_{\sigma}(x) = \min\{\sigma(y) \mid y \in \widehat{Q} \mbox{ covers } x\}$ and $\Delta_{\sigma}(x) = \max\{\sigma(y) \mid y \in \widehat{Q} \mbox{ is covered by } x\}$.
  Define the  \textbf{toggle} $\tau_{x}:\mathcal{A}^{\ell}({Q})\rightarrow \mathcal{A}^{\ell}({Q})$
  by \[ \tau_{x}(\sigma)(x'):= \begin{cases}
      \sigma(x')                                            & x \neq x' \\
      \nabla_{\sigma}(x') + \Delta_{\sigma}(x') -\sigma(x') & x=x'.
    \end{cases} \]
\end{definition}

\begin{remark}
  \label{remark:commute}
  It is easy to see, and well--known, the $\tau_{x}$ satisfy:
  \begin{enumerate}
    \item $\tau_{x}^2=1$, and
    \item $\tau_{x}$ and $\tau_{x'}$ commute whenever  $x$ and $x'$ do not share a covering relation.
  \end{enumerate}
\end{remark}

\begin{definition}
  \label{def:row1}
  \textbf{Rowmotion} on $\mathcal{A}^{\ell}({Q})$ is defined as the toggle composition $\row := \tau_{x_1}\circ \tau_{x_2}\circ \cdots\circ \tau_{x_m}$ where $x_1,x_2,\ldots,x_m$ is any linear extension of $Q$.
\end{definition}

\subsection{$P$-strict promotion and $\Gamma(P,R_{\alpha}^{\beta})$ rowmotion}
\label{sec:results}
In this subsection, we give the main results of \cite{BSV2021} in slightly less generality (the case $u=v=0$), which is all we need for this paper. Since the statements do not match exactly the statements in \cite{BSV2021}, we give the specific references of each statement from \cite{BSV2021} for comparison. In the next subsection, we further specialize these results to the case where $P$ is graded.

We will need the following definition, adapted from \cite[Definition 2.1]{DSV2019}.
\begin{definition}
  A restriction function $R_\alpha^\beta$ is \textbf{consistent} if, for every covering relation $x \lessdot y$ in $P$, we have $\alpha(x) < \alpha(y)$ and $\beta(x) < \beta(y)$. 
\end{definition}

Given a generic restriction function $R_{\alpha}^{\beta}$, let $A,B$ be such that $R_A^B$ is a consistent restriction function with $\alpha(p)\leq A(p)\leq B(p)\leq \beta(p)$ for each $p\in P$ and each $[A(p),B(p)]$ is the largest possible interval. We say $R_A^B$ is the consistent restriction function associated to $R_{\alpha}^{\beta}$.

\begin{definition}[\protect{\cite[Theorem 2.21]{DSV2019}}]
  \label{def:GammaOne}
Given a restriction function $R_{\alpha}^{\beta}$, let $R_A^B$ be the consistent restriction function associated to $R_{\alpha}^{\beta}$. The \textbf{gamma poset} $\Gamma(P,R_{\alpha}^{\beta})$ is  the poset whose elements are $(p,k)$ with $p\in P$ and $A(p)\leq k\leq B(p)-1$ and whose covering relations are given by $(p_1,k_1)\lessdot (p_2,k_2)$ if and only if either
\begin{enumerate}
\item  $p_1=p_2$ and $k_1=k_2+1$, or
\item $p_1\lessdot_P p_2$ and $k_1+1=k_2$.
  \end{enumerate}
\end{definition}

\begin{figure}[hbtp]
  \begin{center}
    \includegraphics[width=\textwidth]{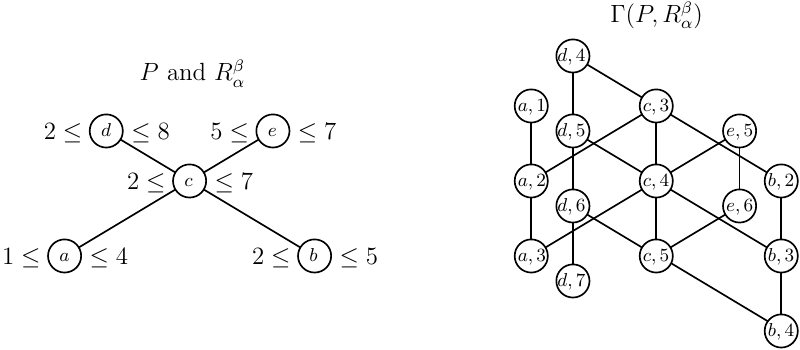}
  \end{center}
  \caption{An example of a gamma poset.  The poset $P$ is drawn on the left along with its restriction function $R_\alpha^\beta$ indicated by the numerical bounds on each element, and the poset $\Gamma(P,R_\alpha^\beta)$ constructed using Definition~\ref{def:GammaOne} is shown on the right.  Observe that $R_\alpha^\beta$ is not consistent, and so $\Gamma(P,R_\alpha^\beta)$ contains only those elements $(p,k)$ where $k \in R_A^B(p)$.  For example, since the least possible label for $b$ is 2 and the greatest possible label for $e$ is 7, we have $A(c) = 3$ and $B(c) = 6$. Thus, $\Gamma(P,R_\alpha^\beta)$ contains the elements $(c,3), (c,4)$ and $(c,5)$.} \label{fig:GammaExX}
\end{figure}

Figure~\ref{fig:GammaExX} shows an example of the gamma poset for the X-shaped poset given a restriction function $R_\alpha^\beta$ which is not consistent.

\begin{remark}
Note that in the above, $R_A^B$ is a consistent restriction function of intervals $[A(p),B(p)]$ on each poset element.  Therefore, we take \cite[Theorem 2.21]{DSV2019} as our definition of the gamma poset, since we do not need the full generality of \cite[Definition 2.11]{DSV2019}. 
\end{remark}

\begin{definition}[\protect{\cite[Definition~2.6]{BSV2021}}]
  \label{def:togpro}
  \textbf{Toggle-promotion} on $\mathcal{A}^{\ell}({\Gamma}(P,R_{\alpha}^{\beta}))$ is defined as the toggle composition $\togpro := \cdots \circ\tau_3 \circ \tau_{2}\circ \tau_{1}$, where $\tau_{k}$ denotes the composition of all the $\tau_{(p,k)}$ over all $p \in P$, $(p,k) \in\Gamma(P,R_{\alpha}^{\beta})$.
\end{definition}

This composition is well-defined, since the toggles within each  $\tau_k$ commute by Remark~\ref{remark:commute}.

\smallskip
The first main result of \cite{BSV2021}, which we state below as Theorem~\ref{thm:moregeneralpro}, gives an equivariant bijection between $P$-strict labelings and certain $Q$-partitions, via the map of Definition~\ref{def:mainbijection}.

\begin{definition}[\protect{\cite[Definition~2.9]{BSV2021}}]\label{def:mainbijection}
  We define the map $\Phi: \mathcal{L}_{P \times [\ell]}(R_\alpha^\beta) \rightarrow \mathcal{A}^{\ell}(\Gamma(P,R_\alpha^\beta))$ as the composition of two intermediate maps $\phi_2$ and $\phi_3$.  Start with a $P$-strict labeling $f \in \mathcal{L}_{P \times [\ell]}(R_\alpha^\beta)$.
  First, $\phi_2$ sends ${f}$ to the multichain $\mathcal{O}_{\ell} \leq \mathcal{O}_{\ell - 1} \leq \cdots \leq \mathcal{O}_1$ in $J(\Gamma(P,R_\alpha^\beta))$ where, for $1 \leq i \leq \ell$ and $L_i$ the $i$th layer of $P \times [\ell]$, $\phi_2$ sends ${f}(L_i)$ to its associated order ideal $\mathcal{O}_i \in J(\Gamma(P,R_\alpha^\beta))$
where $\mathcal{O}_i$ is generated by the set $\{(p, k) ~ | ~ f(p) = k\}$.
Then, $\phi_3$ maps the above multichain to a $\Gamma(P,R_\alpha^\beta)$-partition $\sigma$ as seen in \cite[p.~11]{Stanley1972}, where $\sigma(p,k) = \#\{i \mid (p,k) \notin \mathcal{O}_i\}$, the number of order ideals not including $(p,k)$.  Let $\Phi = \phi_3 \circ \phi_2$.
\end{definition}

\begin{figure}[hbtp]
  \begin{center}
    \includegraphics[width=\textwidth]{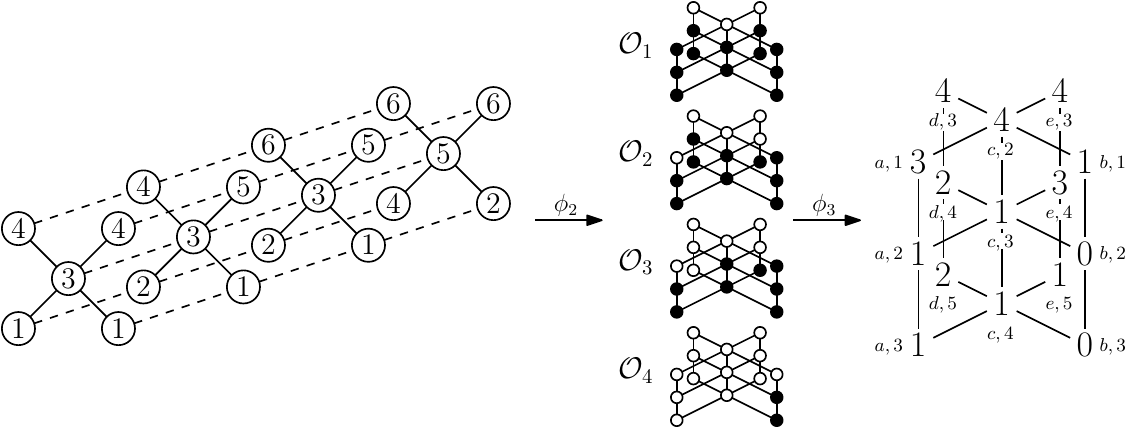}
  \end{center}
  \caption{An example of the map $\Phi$.  On the left is an element $f \in \mathcal{L}_{P \times [\ell]}(R^6)$, where $P$ is the X-shaped poset shown in Figure \ref{fig:GammaExX}.  On the right is $\Phi(f)$, the corresponding element of $\mathcal{A}^{4}(\Gamma(P,R^q))$, which includes the element names $(p,k)$ of $\Gamma(P,R^q)$ for the purpose of clarity.  In the middle is the intermediate result $\phi_2(f)$, the multichain of order ideals $\mathcal{O}_{4} \leq \mathcal{O}_{3} \leq \mathcal{O}_2 \leq \mathcal{O}_1$ in $J(\Gamma(P,R^6))$.  The map $\phi_3$ acts on this multichain to produce the labels of $\Phi(f)$.  For example, the element $(e,4)$ is labeled by a 3 because the three order ideals $\mathcal{O}_2, \mathcal{O}_3,$ and $\mathcal{O}_4$ do not include the element $(e,4)$.} \label{fig:PhiMap}
\end{figure}

Figure \ref{fig:PhiMap} shows an example of the map $\Phi$ via $\phi_2$ and $\phi_3$.

\begin{theorem}[\protect{\cite[Theorem~2.8]{BSV2021}}]
  \label{thm:moregeneralpro}
  $\mathcal{L}_{P\times[\ell]}(R_\alpha^\beta)$ under $\pro$ is in equivariant bijection with $\mathcal{A}^{\ell}(\Gamma(P,R_\alpha^\beta))$ under $\togpro$. More specifically, for $f\in\mathcal{L}_{P\times[\ell]}(R_\alpha^\beta)$, $\Phi\left(\pro(f)\right)=\togpro\left(\Phi(f)\right)$.
\end{theorem}

The equivalence is not only of the actions $\pro$ and $\togpro$, but is proved by corresponding Bender--Knuth involutions with toggles.
\begin{lemma}[\protect{\cite[Lemma~2.11]{BSV2021}}]  \label{lem:equiv}
  The bijection map $\Phi$ equivariantly takes the generalized Bender--Knuth involution $\rho_k$ to the toggle operator $\tau_k$.
\end{lemma}

The second main result of \cite{BSV2021}, Theorem~2.20, states a general condition under which rowmotion is in equivariant bijection with promotion. All the posets we consider in this paper satisfy that condition, thus the following corollary is all that is needed.

\begin{corollary}[\protect{\cite[Corollary~2.24]{BSV2021}}]
  \label{cor:abrow}
  $\mathcal{A}^{\ell}({\Gamma}(P,{R_{\alpha}^{\beta}}))$ under $\row$ is in equivariant bijection with $\mathcal{L}_{P \times [\ell]}(R_{\alpha}^{\beta})$ under $\pro$.
\end{corollary}

\subsection{The graded case: $P$-strict promotion and $P\times[q-n-1]$ rowmotion}
\label{sec:newcor}
In the rest of the paper (with the exception of Section~\ref{sec:resonance}), we concentrate on the case where $P$ is a \textbf{graded} poset of \textbf{rank} $n$, meaning all maximal chains of $P$ have $n+1$ elements. This condition simplifies matters significantly; in particular, it determines $\Gamma(P,R^q)$, the gamma poset with restriction function induced by a global bound $q$. While Lemma~\ref{lem:GradedGlobalq} and Corollary~\ref{cor:PrankedGlobalqRow} were stated and proved in~\cite{BSV2021}, the other results of the subsection are new statements that follow from the general theorems of \cite{BSV2021}. These will be easier to apply to our cases of interest.
\begin{lemma}[\protect{\cite[Lemma~2.28]{BSV2021}}]
  \label{lem:GradedGlobalq}
  Let $P$ be a graded poset of rank $n$.  Then $\Gamma(P,R^q)$ is isomorphic to $P \times [q-n-1]$ as a poset.
\end{lemma}

\begin{remark}
  \label{rem:proof_bij}
  The proof of the above uses the bijection from $\Gamma(P,R^q)$ to $P \times [q-n-1]$ that sends $(p,k)$ to $(p,q-n+\operatorname{rank}(p)-k)$, where the rank of a minimal poset element is zero.
\end{remark}

We now state the versions of our main theorems we will use.

\begin{corollary}
  \label{thm:PrankedGlobalq}
  Let $P$ be a graded poset of rank $n$.  Then $\mathcal{L}_{P \times [\ell]}(R^q)$ under $\pro$ is in equivariant bijection with $\mathcal{A}^\ell(P \times [q-n-1])$ under $\togpro$. More specifically, for $f\in\mathcal{L}_{P\times[\ell]}(R^q)$, $\Phi\left(\pro(f)\right)=\togpro\left(\Phi(f)\right)$.
\end{corollary}
\begin{proof}
  This follows directly from Theorem~\ref{thm:moregeneralpro} and Lemma~\ref{lem:GradedGlobalq}.
\end{proof}

\begin{figure}[hbtp]
  \begin{center}
    \includegraphics[width=\textwidth]{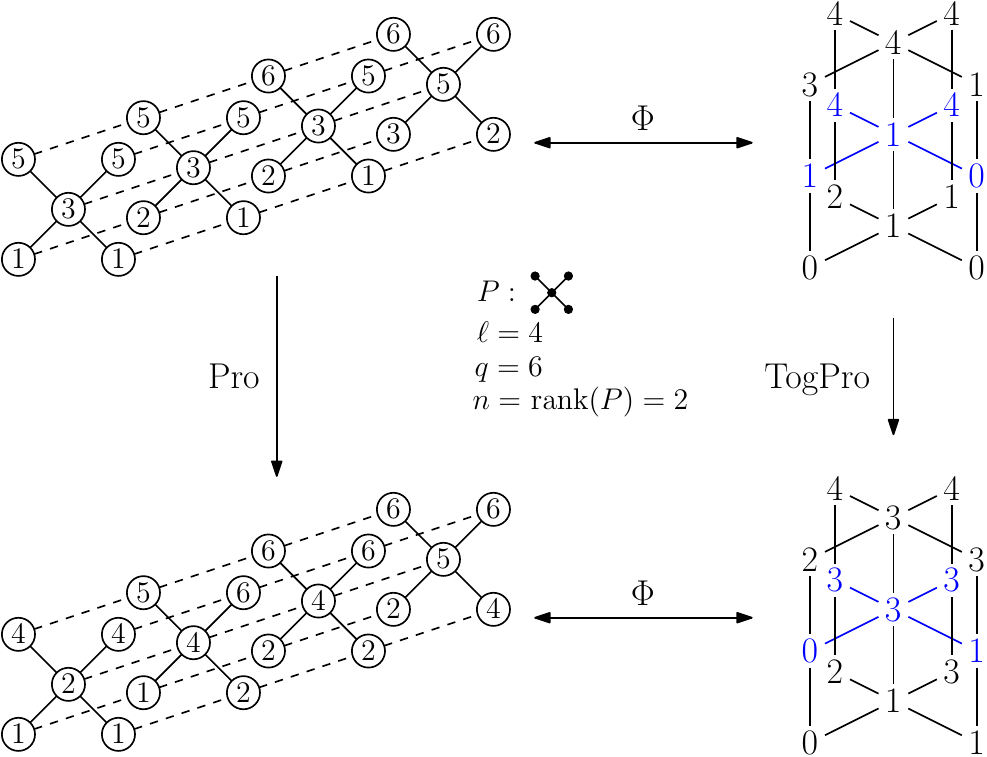}
  \end{center}
  \caption{An illustration of Corollary~\ref{thm:PrankedGlobalq}. Promotion on the $P$-strict labeling in the upper left corresponds to toggle-promotion on the $\Gamma(P,R^q)$-partition in the upper right. The poset $P$ of rank $n$ along with the parameters $\ell$ and $q$ are shown in the center. Since $P$ is graded, $\Gamma(P,R^q)$ is isomorphic to $P \times [q-n-1] = P \times [3]$.  The blue elements of $P \times [3]$ are the elements of the subposet $P \times \{2\}$ described in the statement of Proposition \ref{prop:togpro_rowinverse}.} \label{fig:PrankedGlobalq}
\end{figure}

Figure~\ref{fig:PrankedGlobalq} showcases an example of the bijection from Corollary~\ref{thm:PrankedGlobalq}.

\begin{corollary}[\protect{\cite[Corollary~2.29]{BSV2021}}]
  \label{cor:PrankedGlobalqRow}
  Let $P$ be a graded poset of rank $n$.  Then $\mathcal{L}_{P \times [\ell]}(R^q)$ under $\pro$ is in equivariant bijection with $\mathcal{A}^\ell(P \times [q-n-1])$ under $\row$.
\end{corollary}
\begin{proof}
  This follows directly from Corollary~\ref{cor:abrow} and Lemma~\ref{lem:GradedGlobalq}.
\end{proof}

\begin{remark}
  \label{remark:ell1}
  In the case $\ell=1$, the bijection degenerates to a correspondence from \cite{DSV2019} between order ideals $\mathcal{A}^{1}(P\times[a])=J(P\times [a])$ and $\mathcal{L}_{P\times [1]}(R^{a+n+1})$, which has fibers of size $1$ and is thus equivalent to the set of strictly increasing labelings on $P$ with labels between $1$ and $a+n+1$. When $P$ is graded, this is equivalent (by subtracting $i$ from the labels in rank $i$) to the correspondence of \cite[Prop.\ 3.5.1]{Stanley2011} between order ideals  and order-preserving maps. We will mainly focus on the case $\ell>1$.
\end{remark}

The characterization of the gamma poset in Lemma~\ref{lem:GradedGlobalq} allows us to give an explicit description of the toggles in $\togpro$ without reference to the gamma poset labeling.
\begin{lemma}
  \label{lem:togpro_pbyqn1}
  Let $P$ be a graded poset of rank $n$. Then $\togpro$ on $\Gamma(P,R^q)$ is equivalent to the toggle composition $\togpro := \tau_{q-1} \circ \tau_{q-2} \circ \cdots \circ \tau_{2}\circ \tau_{1}$ on $\mathcal{A}^{\ell}(P \times [q-n-1])$, where $\tau_{k}$ denotes the composition of all the $\tau_{(p,i)}$ over all $p \in P$, where $i = q - n + \operatorname{rank}(p) - k$.
\end{lemma}

\begin{proof}
  Apply the bijection from Remark~\ref{rem:proof_bij} to Definition~\ref{def:togpro}.
\end{proof}

The following gives an equivalent toggle order for $\togpro$, which will be helpful in Corollary~\ref{cor:productofthreechains}. Note that in the proposition below, $\row^{-1}(Q)$ where $Q$ is a subposet of $P$ indicates that we apply the toggles associated to $Q$ in the order of any linear extension of $Q$.
\begin{proposition}
  \label{prop:togpro_rowinverse}
  Let $P$ be a graded poset of rank $n$.  Denote by $P \times \{j\}$ the subposet of $P \times [q-n-1]$ consisting of all elements of the form $(p,j)$, $p \in P$.  Then $\togpro$ on $\Gamma(P,R^q)$ is equivalent to $\row^{-1}(P \times \{1\}) \circ \row^{-1}(P \times \{2\}) \circ \cdots \circ \row^{-1}(P \times \{q-n-2\}) \circ \row^{-1}(P \times \{q-n-1\})$.
\end{proposition}

\begin{proof}
  Let $P_m$ denote the set of elements of $P$ with rank $m$ and, for some $j \in [q-n-1]$, let $\tau_{(P_m,j)}$ be the composition of all the toggles $\tau_{(p,j)}$ where $p \in P_m$.  Then $\row^{-1}(P \times \{j\})$ is given by the toggle composition $\tau_{(P_n,j)} \circ \tau_{(P_{n-1},j)} \circ \cdots \circ \tau_{(P_1,j)} \circ \tau_{(P_0,j)}$, and we have that $\tau_{(P_m,j)}$ commutes with $\tau_{(P_{m'},j')}$ whenever $\vert m-m' \vert + \vert j -j' \vert > 1$.

  Now, the $\tau_k$ in Lemma~\ref{lem:togpro_pbyqn1} are given by \begin{align*}
    \tau_1     & = \tau_{(P_0, q-n-1)}                                                                         \\
    \tau_2     & = \tau_{(P_1, q-n-1)} \circ \tau_{(P_0, q-n-2)}                                               \\
    \tau_3     & = \tau_{(P_2, q-n-1)} \circ \tau_{(P_1, q-n-2)} \circ \tau_{(P_0, q-n-3)}                     \\
               & \shortvdotswithin{=}
    \tau_{n}   & = \tau_{(P_{n-1}, q-n-1)} \circ \tau_{(P_{n-2}, q-n-2)} \circ \cdots \circ \tau_{(P_0, q-2n)} \\
    \tau_{n+1} & = \tau_{(P_n, q-n-1)} \circ \tau_{(P_{n-1}, q-n-2)} \circ \cdots \circ \tau_{(P_0, q-2n-1)}   \\
    \tau_{n+2} & = \tau_{(P_n, q-n-2)} \circ \tau_{(P_{n-1}, q-n-3)} \circ \cdots \circ \tau_{(P_0, q-2n-2)}   \\
    \tau_{n+3} & = \tau_{(P_n, q-n-3)} \circ \tau_{(P_{n-1}, q-n-4)} \circ \cdots \circ \tau_{(P_0, q-2n-3)}   \\
               & \shortvdotswithin{=}
    \tau_{q-3} & = \tau_{(P_n, 3)} \circ \tau_{(P_{n-1}, 2)} \circ \tau_{(P_{n-2}, 1)}                         \\
    \tau_{q-2} & = \tau_{(P_n, 2)} \circ \tau_{(P_{n-1}, 1)}                                                   \\
    \tau_{q-1} & = \tau_{(P_n, 1)},
  \end{align*}
  so we can rewrite $\togpro$ in terms of the $\tau_{(P_m,i)}$.  After doing so, we note that, for all $m>0$, $\tau_{(P_m, q-n-1)}$ commutes with everything to its right except for $\tau_{(P_{m-1}, q-n-1)}$. Therefore, we can reorder $\togpro$ such that all of the $\tau_{(P_m, q-n-1)}$ are on the right as follows:\begin{align*} \togpro = & \tau_{(P_n, 1)} \circ \cdots \circ \tau_{(P_0, q-n-2)} \circ \tau_{(P_n, q-n-1)} \circ \tau_{(P_{n-1}, q-n-1)} \circ \cdots \circ \tau_{(P_1, q-n-1)} \circ \tau_{(P_0, q-n-1)} \\
    =         & \tau_{(P_n, 1)} \circ \cdots \circ \tau_{(P_0, q-n-2)} \circ {\row}^{-1}(P \times \{q-n-1\}).
  \end{align*}

  Similarly, $\tau_{(P_m, q-n-2)}$ commutes with everything to its right except for $\tau_{(P_{m-1}, q-n-2)}$ and ${\row}^{-1}(P \times \{q-n-1\})$, so \[\togpro = \tau_{(P_n, 1)} \circ \cdots \circ \tau_{(P_0, q-n-3)} \circ {\row}^{-1}(P \times \{q-n-2\}) \circ {\row}^{-1}(P \times \{q-n-1\}),\] and so on until we have\[\togpro = {\row}^{-1}(P \times \{1\}) \circ {\row}^{-1}(P \times \{2\}) \circ \cdots \circ {\row}^{-1}(P \times \{q-n-2\}) \circ {\row}^{-1}(P \times \{q-n-1\})\] as desired.
\end{proof}

\section{Products of chains and multifold symmetry}
\label{sec:ab}
Our motivation in this section is to explore the bijection of Corollary~\ref{thm:PrankedGlobalq} further when $P$ is a product of chains. In Subsection~\ref{subsection:applicationofmainresult}, we begin with Corollaries~\ref{cor:productofthreechains} and \ref{cor:productofnchains}, which specialize Corollary~\ref{thm:PrankedGlobalq} when $P$ is a product of chains and states TogPro in terms of the hyperplane toggle definition of \cite{DPS2017}.
In the rest of the section, we use symmetry to prove equivalences of $P$-strict labelings and apply the resulting bijections to obtain order and homomesy results on the $P$-strict labelings $\mathcal{L}_{([a] \times [b]) \times [\ell]}(R^{a+b})$.
Specifically, in Subsection~\ref{subsection:multifoldsymmetry}, we give multifold symmetry results on $P$-strict labelings in Theorems~\ref{thm:3dmultifold} and \ref{thm:ndmultifold}. Because semistandard tableaux can be expressed as certain $P$-strict labelings, we present several results on semistandard Young tableaux in Subsection~\ref{subsection:SSYT}. This includes a new distribution result on semistandard Young tableaux under promotion (Theorem~\ref{ssytbcdist}) and a new homomesy result on semistandard Young tableaux under promotion (Corollary~\ref{cor:boxcounthomomesy}). In Subsection~\ref{subsection:consequencesoftrifold}, we use results on semistandard Young tableaux and our multifold symmetry to establish an order result, a distribution result, and a homomesy result on $\mathcal{L}_{([a] \times [b]) \times [\ell]}(R^{a+b})$ in Corollary~\ref{cor:aborder}, Theorem~\ref{thm:pstrictdist}, and Corollary~\ref{cor:pstrict_homomesy}, respectively. Finally, we state a conjecture in Subsection~\ref{subsec:conj} and translate it using previous results of this section.

\subsection{Application of main result to $P=[a]\times[b]$}
\label{subsection:applicationofmainresult}

Using the hyperplane toggle definition of \cite{DPS2017}, we can determine which hyperplane sweep on the product of chains poset corresponds to TogPro.

\begin{definition}{\cite[Definition 3.13]{DPS2017}}
  We say that an \emph{$n$-dimensional lattice projection} of a ranked poset $P$ is an order and rank preserving map $\pi : P \rightarrow \mathbb{Z}^n$, where the rank function on $\mathbb{Z}^n$ is the sum of the coordinates and $x \le y$ in $\mathbb{Z}^n$ if and only if the componentwise difference $y-x$ is in $(\mathbb{Z}_{\ge 0})^n$.
\end{definition}

\begin{definition}{\cite[Definition 3.14]{DPS2017}}
  Let $Q$ be a poset with an $n$-dimensional lattice projection $\pi$ and let $v  \in \{\pm 1\}^n$.  Let $T_{\pi, v}^i$ be the product of toggles $t_x$ for all elements $x$ of $Q$ that lie on the affine hyperplane $\langle \pi(x),v \rangle=i$.  If there is no such $x$, then this is the empty product, considered to be the identity.  Define \textit{promotion with respect to $\pi$ and $v$} as the (finite) toggle product Pro$_{\pi,v}=\dots T_{\pi,v}^{-2} T_{\pi,v}^{-1} T_{\pi,v}^{0} T_{\pi,v}^{1} T_{\pi,v}^{2}\dots$.
\end{definition}

\begin{remark}
  For a $Q$-partition $\mathcal{A}^\ell([a_1] \times \dots \times [a_k])$, we will frequently use the identity map for $\pi$; the identity map will be denoted with $\mathrm{id}$. For $Q$-partitions obtained from a $\Gamma$ poset construction, the last coordinate decreases as we traverse up the poset. As a result, we use a non-identity map $\pi$ to compensate for this labeling. See Corollaries~\ref{cor:productofthreechains} and \ref{cor:productofnchains} for instances of this.
\end{remark}

\begin{corollary}
  \label{cor:productofthreechains}
  $\mathcal{L}_{([a] \times [b]) \times [\ell]}(R^{a+b+c-1})$ under $\pro$ is in equivariant bijection with $\mathcal{A}^\ell([a] \times [b] \times [c])$ under $\row$. More specifically, for $f\in\mathcal{L}_{([a] \times [b]) \times [\ell]}(R^{a+b+c-1})$, $\Phi\left(\pro(f)\right)=\mathrm{Pro}_{\pi,(-1,-1,1)}\left(\Phi(f)\right)$ where $\pi((i,j),k)=(i,j,i+j-k+c-1)$.
\end{corollary}

\begin{proof}
  By setting $P=[a] \times [b]$ and $q=a+b+c-1$, the existence of an equivariant bijection follows directly from Corollary~\ref{thm:PrankedGlobalq}. To show that $\Phi\left(\pro(f)\right)=\mathrm{Pro}_{\pi,(-1,-1,1)}\left(\Phi(f)\right)$, we will show that $\mathrm{TogPro}$ and $\mathrm{Pro}_{\pi,(-1,-1,1)}$ have the same toggle order on $\Gamma(P,R^q)$. By Theorem 23 of \cite{Vorland2019}, the toggles of $\mathrm{Pro}_{\pi,(-1,-1,1)}$ can be reordered to $\row^{-1}([a] \times [b] \times \{1\}) \circ \row^{-1}([a] \times [b] \times \{2\}) \circ \cdots \circ \row^{-1}([a] \times [b] \times \{c-1\}) \circ \row^{-1}([a] \times [b] \times \{c\})$. By Proposition~\ref{prop:togpro_rowinverse}, this is equivalent to $\togpro$. It follows that $\Phi\left(\pro(f)\right)=\mathrm{Pro}_{\pi,(-1,-1,1)}\left(\Phi(f)\right)$ from Corollary~\ref{thm:PrankedGlobalq}.
\end{proof}

\begin{remark}
  In \cite{DPS2017}, K.\ Dilks, O.\ Pechenik, and Striker showed that if $I$ is an order ideal of $[a] \times [b]\times [c]$, the action $\mathrm{Pro}_{(1,1,-1)}$ on $I$ corresponds to $K$-promotion on a corresponding increasing tableaux. The reason $\mathrm{Pro}_{(1,1,-1)}$ is used in their result is because of labeling conventions: an element in an order ideal is labeled 1 and an element not in an order ideal is labeled 0. With $Q$-partitions, we are using the reverse convention, which is the reason for $\mathrm{Pro}_{\pi,(-1,-1,1)}$ appearing in Corollary~\ref{cor:productofthreechains}.
\end{remark}

The following corollary extends Corollary~\ref{cor:productofthreechains} to products of more than two chains. Its proof is nearly identical to the proof of Corollary~\ref{cor:productofthreechains}, so we omit it here. By setting $q=\sum_i^k a_i-1$, the existence of an equivariant bijection follows directly from Corollary~\ref{thm:PrankedGlobalq}.
\begin{corollary}
  \label{cor:productofnchains}
  Let $P=[a_1] \times [a_2] \times \dots \times [a_{k-1}]$. Then $\mathcal{L}_{P \times [\ell]}(R^{\sum_i^k a_i-1})$ under $\pro$ is in equivariant bijection with $\mathcal{A}^\ell(P \times [a_k])$ under $\row$. More specifically, for $f\in\mathcal{L}_{P \times [\ell]}(R^{\sum_i^k a_i-1})$, $\Phi\left(\pro(f)\right)=\mathrm{Pro}_{(-1,\dots, -1, 1)}\left(\Phi(f)\right)$ where $\pi((i_1,i_2,\dots,i_{k-1}),i_k)=(i_1,i_2,\dots,i_{k-1}, i_1+\dots+i_{k-1}-i_k+a_k-1)$.
\end{corollary}

\subsection{Multifold symmetry}
\label{subsection:multifoldsymmetry}
In this subsection, we use the product of chains results of Corollaries~\ref{cor:productofthreechains} and~\ref{cor:productofnchains} to show $P$-strict trifold symmetry and $P$-strict multifold symmetry in Theorems~\ref{thm:3dmultifold} and \ref{thm:ndmultifold}, respectively. We also state a similar, more general equivalence in Theorem~\ref{thm:twofoldsymmetry}.

\begin{theorem}
  \label{thm:3dmultifold}
  There are promotion-equivariant bijections among the sets $\mathcal{L}_{([a] \times [b]) \times [\ell]}(R^{a+b+c-1})$, $\mathcal{L}_{([a] \times [c]) \times [\ell]}(R^{a+b+c-1})$, and $\mathcal{L}_{([b] \times [c]) \times [\ell]}(R^{a+b+c-1})$.
\end{theorem}

\begin{proof}
  By Corollary~\ref{cor:productofthreechains}, $\mathcal{L}_{([a] \times [b]) \times [\ell]}(R^{a+b+c-1})$ under $\pro$ is in equivariant bijection with $\mathcal{A}^\ell([a] \times [b] \times [c])$ under $\row$, $\mathcal{L}_{([a] \times [c]) \times [\ell]}(R^{a+b+c-1})$ under $\pro$ is in equivariant bijection with $\mathcal{A}^\ell([a] \times [c] \times [b])$ under $\row$, and $\mathcal{L}_{([b] \times [c]) \times [\ell]}(R^{a+b+c-1})$ under $\pro$ is in equivariant bijection with $\mathcal{A}^\ell([b] \times [c] \times [a])$ under $\row$. The set of linear extensions of $[a]\times[b]\times[c]$ is the same under any permutation of $a,b,c$. As a result, all three sets under their respective promotion actions are in equivariant bijection with $\mathcal{A}^\ell([a] \times [b] \times [c])$ under $\row$, and so there are equivariant bijections between the three sets under their respective promotion actions.
\end{proof}

We can extend this same idea to a more general product of chains in the theorem below. Since the proof is nearly identical to the proof of Theorem~\ref{thm:3dmultifold}, we give it in condensed form.

\begin{theorem}
  \label{thm:ndmultifold}
  Let $P_i=[a_1] \times [a_2] \times \dots \times [a_{i-1}] \times [a_{i+1}] \times \dots \times [a_{k}]$. For any $1 \le i,j \le k$, there is an equivariant bijections between the sets $\mathcal{L}_{P_i \times [\ell]}(R^{\sum_i^k a_i-1})$ and $\mathcal{L}_{P_j \times [\ell]}(R^{\sum_i^k a_i-1})$ under their respective promotion actions.
\end{theorem}
\begin{proof}
  Both sets are in equivariant bijection with $\mathcal{A}^\ell([a_1] \times \dots \times [a_k])$ under Row by Corollary~\ref{cor:productofnchains}.
\end{proof}

In the case of the product of an arbitrary graded poset with a chain, we have a similar symmetry, but the equivalence is between only two sets rather than three or more.

\begin{theorem}
  \label{thm:twofoldsymmetry}
  Let $P$ be graded poset of rank $n$. Then there is an equivariant bijection between $\mathcal{L}_{(P \times [a]) \times [\ell]}(R^{a+b+n})$ under $\pro$ and $\mathcal{L}_{(P \times [b]) \times [\ell]}(R^{a+b+n})$ under $\pro$.
\end{theorem}

\begin{proof}
  By Corollary~\ref{thm:PrankedGlobalq}, both sets are in equivariant bijection with $\mathcal{A}^\ell([P] \times [a] \times [b])$ under $\row$.
\end{proof}

\begin{remark}
  When $a=1$, $R^{1+b+n}$ restricts the labels of $\mathcal{L}_{(P \times [b]) \times [\ell]}(R^{1+b+n})$ to only two possible values at any $(p,i)$.  Thus, we can think of the layers of $f \in \mathcal{L}_{(P \times [b]) \times [\ell]}(R^{1+b+n})$ as order ideals of $P \times [b]$ by considering all elements labeled by the lower of their two values as elements of the corresponding order ideal.  Then, $f$ itself corresponds to a multichain of order ideals of $P \times [b]$, and we can interpret the $a=1$ case of Theorem~\ref{thm:twofoldsymmetry} as an application of the intermediate bijection $\phi_2$ from Definition~\ref{def:mainbijection}.
\end{remark}

\subsection{Semistandard Young Tableaux}
\label{subsection:SSYT}

In this subsection, we specialize the trifold symmetry of Theorem~\ref{thm:3dmultifold} to certain $P$-strict labelings and semistandard Young tableaux. We also
obtain new distribution and homomesy results on semistandard Young tableaux (Theorem~\ref{ssytbcdist} and Corollary~\ref{cor:boxcounthomomesy}). In the next subsection, we use Theorem~\ref{ssytbcdist} in conjunction with our main bijection to prove Theorem~\ref{thm:ppartdist}, a distribution result on $Q$-partitions on a product of two chains.

\begin{definition}
  Let $\mathrm{SSYT}_{k}(m \times n)$ denote the set of semistandard Young tableaux of shape $m \times n$ with largest entry $k$.
\end{definition}

\begin{corollary}
  \label{cor:pstrict_ssyt_bij}
  There are equivariant bijections between $\mathcal{L}_{([a] \times [b]) \times [\ell]}(R^{a+b})$, $\mathrm{SSYT}_{a+b}(a \times \ell)$, and $\mathrm{SSYT}_{a+b}(b \times \ell)$ under their respective promotion actions.
\end{corollary}

\begin{proof}
  By specializing to $c=1$ into Theorem \ref{thm:3dmultifold}, we obtain equivariant bijections between the sets $\mathcal{L}_{([a] \times [b]) \times [\ell]}(R^{a+b})$, $\mathcal{L}_{([a] \times [1]) \times [\ell]}(R^{a+b})$, and $\mathcal{L}_{([b] \times [1]) \times [\ell]}(R^{a+b})$ under their respective promotion actions. The set $\mathcal{L}_{([a] \times [1]) \times [\ell]}(R^{a+b})$ is in bijection with $\mathrm{SSYT}_{a+b}(a \times \ell)$ and the set $\mathcal{L}_{([b] \times [1]) \times [\ell]}(R^{a+b})$ is in bijection with $\mathrm{SSYT}_{a+b}(b \times \ell)$, completing the proof.
\end{proof}

Corollary \ref{cor:pstrict_ssyt_bij} shows there is a bijection between $\mathrm{SSYT}_{a+b}(a \times \ell)$ and $\mathcal{L}_{([a] \times [b]) \times [\ell]}(R^{a+b})$. From \cite{BPS2016}, there are known distribution and homomesy results on semistandard Young tableaux under promotion. We begin with the definition of homomesy and a relevant result of J.\ Bloom, Pechenik, and D.\ Saracino.

\begin{definition}[\cite{PR2015}]
  \label{def:homomesy}
  Given a finite set $S$, an action $\tau:S \rightarrow S$, and a statistic $f:S \rightarrow \mathbb{C}$, we say that $(S, \tau, f)$ exhibits \textbf{homomesy} if there exists $c \in \mathbb{C}$ such that for every $\tau$-orbit $\orb$
  \begin{center}
    $\displaystyle\frac{1}{|\orb|} \sum_{x \in \orb} f(x) = c$
  \end{center}
  where $|\orb|$ denotes the number of elements in $\orb$. If such a $c$ exists, we say the triple is \textbf{$c$-mesic}.
\end{definition}

\begin{theorem}[\protect{\cite[Theorem 1.1]{BPS2016}}]
  \label{thm:bpshomomesy}
  Let $S$ denote a set of boxes in an $a \times b$ rectangle fixed under $180^\circ$ rotation, and let $\chi_S(T)$ denote the sum of the entries in the boxes of $S$ in tableaux $T$. Then $(\mathrm{SSYT}_k(a \times b),\mathrm{Pro}, \chi_S)$ is homomesic.
\end{theorem}

To prove Theorem~\ref{thm:bpshomomesy}, Bloom, Pechenik, and Saracino used the fact that the order of $\pro$ on $\mathrm{SSYT}_{k}(a \times b)$ is $k$ (see e.g.~\cite[Theorem 2.9]{BPS2016}) and the notion of distribution on a semistandard Young tableaux under promotion. We give the definition of the distribution and their relevant result below.

\begin{definition}[\protect{\cite[Definition 3.1]{BPS2016}}]
  \label{def:ssytdist}
  Let $T \in \mathrm{SSYT}_k(a \times b)$. For a box $B$ in $T$, define $\mathrm{Dist}(T,B)$ to be the multiset
  \[\mathrm{Dist}(T,B)= \{\chi_{\{B\}}(\mathrm{Pro}^i(T)) : 0 \le i \le k-1 \}. \]
\end{definition}

\begin{proposition}[\protect{\cite[Proof of Theorem 1.1]{BPS2016}}]
  \label{prop:ssytdist}
  Let $T \in \mathrm{SSYT}_k(a \times b)$, $B$ a box in $T$, and $B^*$ the box corresponding to $B$ under $180^\circ$ rotation. Then $\mathrm{Dist}(T,B)=\{k+1-m : m \in \mathrm{Dist}(T,B^*) \}$.
\end{proposition}

Using Proposition \ref{prop:ssytdist}, we prove a new distribution result on semistandard Young tableaux.

\begin{definition}
  Let $T \in \mathrm{SSYT}_k(a \times b)$. For a row $x$ in $T$ and $d \in \mathbb{Z}$, define the box count distribution $\mathrm{BCDist}(T,x,d)$ to be the multiset
  \[\mathrm{BCDist}(T,x,d)= \{ \# \{ (x,j) \in a \times b :  \mathrm{Pro}^i(T)(x,j) > d \} : 0 \le i \le k-1 \}. \]
\end{definition}

The following theorem is a box count distribution result on semistandard Young tableaux under promotion. This result yields a new homomesy result on semistandard Young tableaux under promotion (Corollary~\ref{cor:boxcounthomomesy}), and is used in the proof of a distribution on $Q$-partitions under $\pro_{\mathrm{id},v}$ (Theorem~\ref{thm:ppartdist}).

\begin{theorem}
  \label{ssytbcdist}
  Let $T \in \mathrm{SSYT}_k(a \times b)$, $x$ a row in $T$, and $d \in \mathbb{Z}$. Then $\mathrm{BCDist}(T,a+1-x,k-d)=\{b-m : m \in \mathrm{BCDist}(T,x,d) \}$.
\end{theorem}

\begin{proof}
  We consider each $m \in \mathrm{BCDist}(T,x,d)$ for $m \in \{0, \dots, M\}$, where $M$ denotes the largest value that appears in $\mathrm{BCDist}(T,x,d)$. We split the proof into three cases.

  \emph{Case: $m \ne 0$ and $m \ne b$. } Suppose there are $n_m$ values of $m$ in $\mathrm{BCDist}(T,x,d)$. For each of these, there is a unique $i$ such that $m=\# \{ (x,j) \in a \times b :  \mathrm{Pro}^i(T)(x,j) > d \}$. This implies that for these $n_m$ values of $i$, $\mathrm{Pro}^i(T)(x,b+1-m) > d$ and $\mathrm{Pro}^i(T)(x,b-m) \le d$. Note that there are $n_{m+1}+\dots+n_{M}$ additional $i$ such that $\mathrm{Pro}^{i}(T)(x,b+1-m) > d$ and $\mathrm{Pro}^{i}(T)(x,b-m) > d$. These are the tableaux that have more than $m$ entries larger than $d$ in row $x$. As a result, we have $n_m + n_{m+1}+\dots+n_{M}$ unique $i$ that satisfy $\mathrm{Pro}^{i}(T)(x,b+1-m) > d$. By Proposition~\ref{prop:ssytdist}, this implies we have exactly $n_m + n_{m+1}+\dots+n_{M}$ unique $i'$ that satisfy $\mathrm{Pro}^{i'}(T)(a+1-x,m) < k + 1 - d$. Applying Proposition~\ref{prop:ssytdist} again, we must have exactly $n_m$ values of $i'$ such that $\mathrm{Pro}^{i'}(T)(a+1-x,m+1) \ge k + 1 - d$ and exactly $n_{m+1}+\dots+n_{M}$ values of $i'$ such that $\mathrm{Pro}^{i'}(T)(a+1-x,m+1) < k + 1 - d$. Because there are exactly $n_m$ values of $i'$ such that $\mathrm{Pro}^{i'}(T)(a+1-x,m) < k + 1 - d$ and $\mathrm{Pro}^{i'}(T)(a+1-x,m+1) \ge k + 1 - d$, we have $n_m$ values of $i'$ such that $\mathrm{Pro}^{i'}(T)(a+1-x,m) \le k - d$ and $\mathrm{Pro}^{i'}(T)(a+1-x,m+1) > k - d$. This implies there are $n_m$ values of $b-m$ in $\mathrm{BCDist}(T,a+1-x,k-d)$.

  \emph{Case: $m = 0$. } Suppose there are $n_0$ values of $0$ in $\mathrm{BCDist}(T,x,d)$. For each of these, there is a unique $i$ such that $0=\# \{ (x,j) \in a \times b :  \mathrm{Pro}^i(T)(x,j) > d \}$. This implies that for these $n_0$ values of $i$, $\mathrm{Pro}^i(T)(x,b) \le d$. Note that for all other $i$, we have $\mathrm{Pro}^{i}(T)(x,b) > d$. Applying Proposition~\ref{prop:ssytdist} to box $(x,b)$, we must have exactly $n_0$ values of $i'$ such that $\mathrm{Pro}^{i'}(T)(a+1-x,1) \ge k + 1 - d$, and for all other $i'$, we have $\mathrm{Pro}^{i'}(T)(a+1-x,1) < k + 1 - d$. Therefore, we have exactly $n_0$ values of $i'$ such that $\mathrm{Pro}^{i'}(T)(a+1-x,1) > k - d$. This implies there are $n_0$ values of $b$ in $\mathrm{BCDist}(T,a+1-x,k-d)$.

  \emph{Case: $m=b$. } Suppose there are $n_b$ values of $b$ in $\mathrm{BCDist}(T,x,d)$. For each of these, there is a unique $i$ such that $b=\# \{ (x,j) \in a \times b :  \mathrm{Pro}^i(T)(x,j) > d \}$. This implies that for these $n_b$ values of $i$, $\mathrm{Pro}^i(T)(x,1) > d$. Note that for all other $i$, we have, $\mathrm{Pro}^i(T)(x,1) \le d$. Applying Proposition~\ref{prop:ssytdist} to box $(x,1)$, we must have exactly $n_b$ values of $i'$ such that $\mathrm{Pro}^{i'}(T)(a+1-x,b) < k + 1 - d$, and for all other $i'$, we have $\mathrm{Pro}^{i'}(T)(a+1-x,b) \ge k + 1 - d$. Therefore, we have exactly $n_b$ values of $i'$ such that $\mathrm{Pro}^{i'}(T)(a+1-x,b) \le k - d$. This implies there are $n_b$ values of $0$ in $\mathrm{BCDist}(T,a+1-x,k-d)$.

  Because the theorem statement holds for each $m \in \{0, \dots, M\}$, this completes the proof.
\end{proof}

\begin{figure}[htbp]
  \begin{center}
    \includegraphics[width=\textwidth]{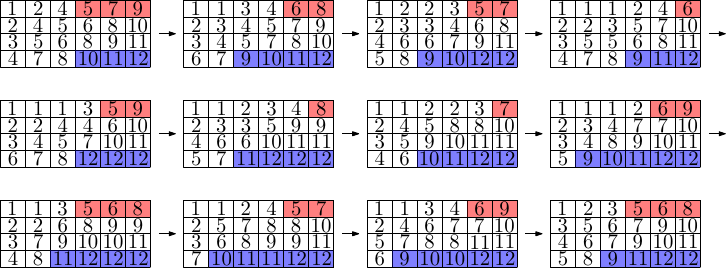}
  \end{center}
  \caption{This figure refers to Example~\ref{ex:boxcount}, which demonstrates Theorem~\ref{ssytbcdist}. From the red boxes, we have $\mathrm{BCDist}(T,1,4)=\{1, 1, 1, 2, 2, 2, 2, 2, 2, 3, 3, 3\}$ while from the blue boxes, we have $\mathrm{BCDist}(T,4,8)=\{3, 3, 3, 4, 4, 4, 4, 4, 4, 5, 5, 5\}$. A quick check verifies $\{6-m : m \in \mathrm{BCDist}(T,1,4) \}=\{3, 3, 3, 4, 4, 4, 4, 4, 4, 5, 5, 5\}=\mathrm{BCDist}(T,4,8)$.}
  \label{fig:boxcountexample}
\end{figure}

\begin{example}
  \label{ex:boxcount}
  Figure~\ref{fig:boxcountexample} is an example of Theorem~\ref{ssytbcdist} where $T \in \mathrm{SSYT}_{12}(4 \times 6)$, $x=1$, and $d=4$. In red, we colored all boxes that would contribute towards the cardinality $\# \{ (1,j) \in 4 \times 6 :  \mathrm{Pro}^i(T)(1,j) > 4 \}$. Recording these cardinalities gives us \[\mathrm{BCDist}(T,1,4)=\{3, 2, 2, 1, 2, 1, 1, 2, 3, 2, 2, 3\}=\{1, 1, 1, 2, 2, 2, 2, 2, 2, 3, 3, 3\}.\] To determine $\mathrm{BCDist}(T,4,8)$, we record the cardinalities of the boxes colored blue, obtaining \[\mathrm{BCDist}(T,4,8)=\{3, 4, 4, 3, 3, 4, 4, 5, 4, 5, 5, 4\}=\{3, 3, 3, 4, 4, 4, 4, 4, 4, 5, 5, 5\}.\] Comparing $\mathrm{BCDist}(T,4,8)$ and $\{6-m : m \in \mathrm{BCDist}(T,1,4) \}=\{3, 3, 3, 4, 4, 4, 4, 4, 4, 5, 5, 5\}$ reveals that these are the same, demonstrating the theorem.
\end{example}

We use Theorem~\ref{ssytbcdist} to prove a new homomesy result on semistandard Young tableaux under promotion. The homomesic statistic is a sum of two simpler statistics; this arises from the two rows in the $\mathrm{BCDist}$ result of Theorem~\ref{ssytbcdist}.

\begin{definition}
  Let $T \in \mathrm{SSYT}_{k}(a \times b)$. Define the statistic $\#B(T,x,d)$ to be the number of boxes in row $x$ of $T$ with entry larger than $d$. In other words,
  \[\#B(T,x,d) = \# \{ (x,j) \in a \times b :  T(x,j) > d \}.\]
\end{definition}

\begin{corollary}
  \label{cor:boxcounthomomesy}
  For any $x \in \{1,\dots,a\}$ and any $d \in \mathbb{Z}$, \[\left(\mathrm{SSYT}_{k}(a \times b), \pro, \#B(T,x,d)+ \#B(T,a+1-x,k-d)\right)\] is $b$-mesic.
\end{corollary}

\begin{proof}
  From Theorem~2.9 of \cite{BPS2016}, the order of $\pro$ on $\mathrm{SSYT}_{k}(a \times b)$ is $k$. Because of this, we can determine orbit averages of $\#B(T,x,d)+ \#B(T,a+1-x,k-d)$ by averaging this statistic over $k$ applications of $\pro$. Taking this approach in conjunction with Theorem~\ref{ssytbcdist} gives us
  \begin{align*}
    \frac{1}{k}\sum_{0 \le i \le k-1}\left(\#B(\mathrm{Pro}^{i}(T),x,d)+\#B(\mathrm{Pro}^{i}(T),a+1-x,k-d) \right) & =  \frac{1}{k} \sum_{m \in \mathrm{BCDist}(T,x,d)}  \left(m + b-m \right) \\ &= \frac{bk}{k} = b.
  \end{align*}
\end{proof}

\subsection{Consequences of multifold symmetry: order and homomesy}
\label{subsection:consequencesoftrifold}
We now use multifold symmetry in conjunction with results on semistandard Young tableaux to obtain results on certain $P$-strict labelings and $Q$-partitions. We begin by applying Corollary \ref{cor:pstrict_ssyt_bij} to obtain the order of $\mathcal{L}_{([a] \times [b]) \times [\ell]}(R^{a+b})$ under promotion. In Definitions~\ref{def:pstrictdist} and \ref{def:ppartdist}, we adapt Definition~\ref{def:ssytdist} on the distribution of a semistandard Young tableaux to the $P$-strict labeling  and $Q$-partition settings, respectively. Using our result on the box count distribution of semistandard Young tableaux under promotion (Theorem~\ref{ssytbcdist}), we are able to give a new proof to a known distribution result on $Q$-partitions under $\pro_{\pi, v}$ (Theorem~\ref{thm:ppartdist}) and prove a new distribution result on $P$-strict labelings under $\pro$ (Theorem~\ref{thm:pstrictdist}). The latter theorem immediately yields Corollary~\ref{cor:pstrict_homomesy}, a homomesy result on $\mathcal{L}_{([a] \times [b]) \times [\ell]}(R^{a+b})$ under promotion with a label sum statistic on antipodal poset elements.

\begin{corollary}
  \label{cor:aborder}
  The order of $\mathcal{L}_{([a] \times [b]) \times [\ell]}(R^{a+b})$ under $\mathrm{Pro}$ divides $a+b$.
\end{corollary}

\begin{proof}
  As shown in Theorem~2.9 of \cite{BPS2016}, the order of $\mathrm{SSYT}_{a+b}(a \times \ell)$ under promotion divides $a+b$. Using Corollary~\ref{cor:pstrict_ssyt_bij}, the order of $\mathcal{L}_{([a] \times [b]) \times [\ell]}(R^{a+b})$ under $\mathrm{Pro}$ divides $a+b$.
\end{proof}

\begin{figure}[htbp]
  \begin{center}
    \includegraphics[width=.90\textwidth]{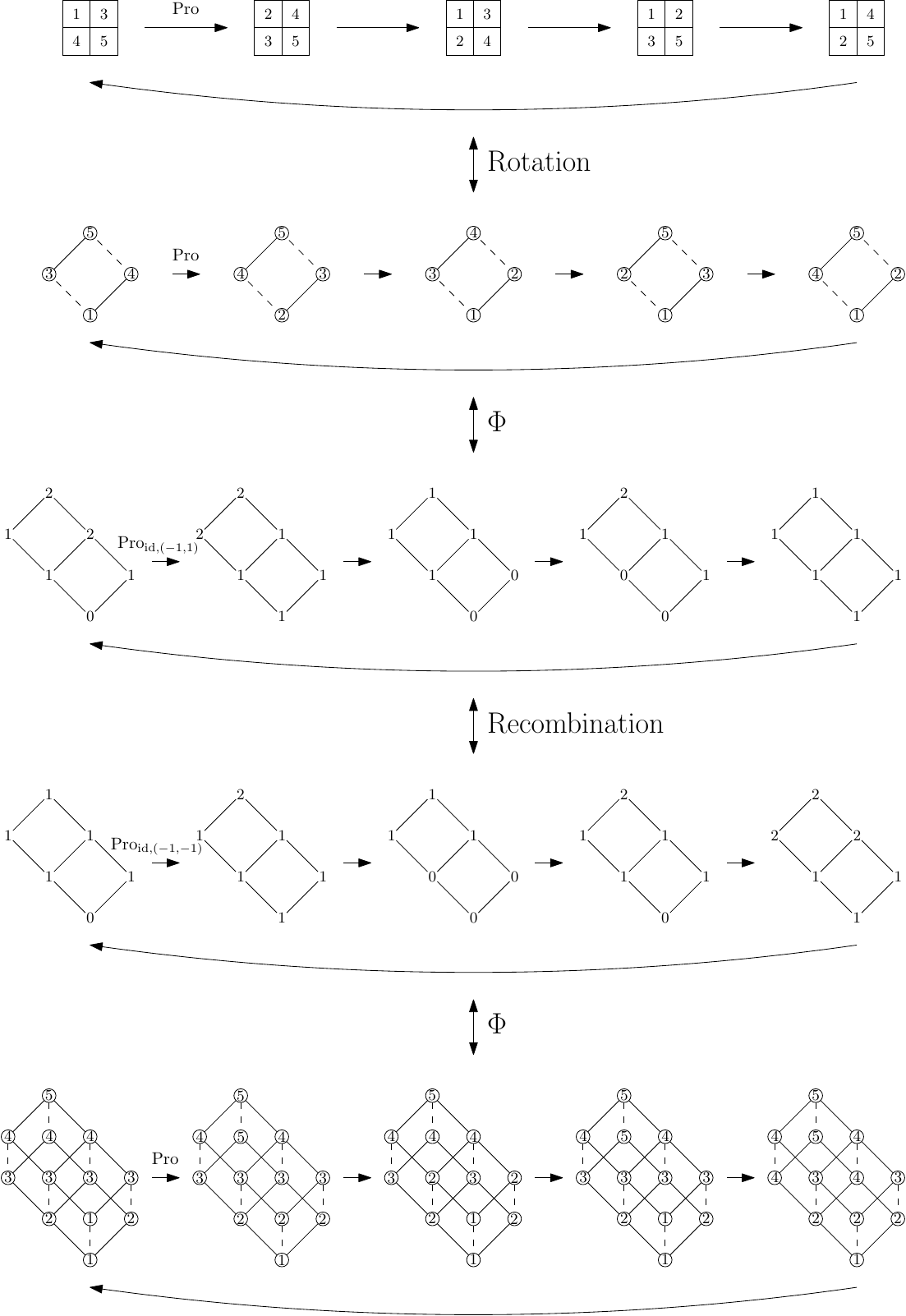}
  \end{center}
  \caption{This figure shows how we can use $\Phi$ to create a bijection from semistandard Young tableaux under $\pro$ to $\mathcal{L}_{([a] \times [b]) \times [\ell]}(R^{a+b})$ under $\pro$. We use this in Theorems~\ref{thm:ppartdist} and \ref{thm:pstrictdist} to translate the distribution result of Theorem~\ref{ssytbcdist} on $\mathrm{SSYT}_{a+c}(a \times \ell)$ to the distribution result of Theorem~\ref{thm:pstrictdist} on $\mathcal{L}_{([a] \times [b]) \times [\ell]}(R^{a+b})$. See Examples~\ref{ex:orbitlength}, \ref{ex:sstytoppart}, and \ref{ex:pparttopstrict} for additional details.}
  \label{fig:pstricthomomesyexample}
\end{figure}

\begin{example}
  \label{ex:orbitlength}
  The top line of Figure~\ref{fig:pstricthomomesyexample} shows an example orbit from $\mathrm{SSYT}_{5}(2 \times 2)$ under promotion. It also shows a bijection to an orbit of $\mathcal{L}_{([2] \times [3]) \times [2]}(R^{5})$ under $\mathrm{Pro}$, which is the bottom line of Figure~\ref{fig:pstricthomomesyexample}. Observe that the order for both orbits is $a+b=5$. We will explain Figure~\ref{fig:pstricthomomesyexample} further in Examples~\ref{ex:sstytoppart} and \ref{ex:pparttopstrict}.
\end{example}

\begin{remark}
  A next natural step would be to specialize Theorem~\ref{thm:3dmultifold} to $c=2$. However, we find that even in the case $a=b=c=\ell=2$, rowmotion has 30 orbits of size 5 and 2 orbits of size 9, indicating an unpleasant orbit structure.
\end{remark}

To state our desired homomesy result on $\mathcal{L}_{([a] \times [b]) \times [\ell]}(R^{a+b})$, we first define a notion of antipodal elements in a product of chains poset and give analogous definitions of distribution for $P$-strict labelings and $Q$-partitions.

\begin{definition}
  \label{def:antipode}
  Given $(i_1,i_2,\ldots,i_k)\in [a_1] \times [a_2]\times \cdots \times [a_k]$, let $(a_1+1-i_1,a_2+1-i_2,\ldots,a_k+1-i_k)$ be the  \textbf{antipode} of $(i_1,i_2,\ldots,i_k)$.
\end{definition}

As in the case of semistandard tableaux, we use $\chi_S(f)$ (resp.~$\chi_S(\sigma)$) to denote the sum of labels of a subset of elements $S$ in a $P$-strict labeling $f$ (resp.~$Q$-partition $\sigma$).

\begin{definition}
  \label{def:pstrictdist}
  Let $f \in \mathcal{L}_{([a] \times [b]) \times [\ell]}(R^{a+b+c-1})$. For a poset element $p \in [a] \times [b] \times [\ell]$, define the distribution $\mathrm{Dist}(f,p)$ to be the multiset
  \[\mathrm{Dist}(f,p)= \{\chi_{\{p\}}(\mathrm{Pro}^i(f)) : 0 \le i \le a+b-1 \}. \]
\end{definition}

\begin{definition}
  \label{def:ppartdist}
  Let $\sigma \in \mathcal{A}^\ell([a] \times [b])$. For a poset element $p \in [a] \times [b]$, define the distribution $\mathrm{Dist}(\sigma,p,\pi,v)$ to be the multiset
  \[\mathrm{Dist}(\sigma,p,\pi,v)= \{\chi_{\{p\}}(\mathrm{Pro}_{\pi,v}^i(\sigma)) : 0 \le i \le a+b-1 \}. \]
\end{definition}

We now state Theorem~\ref{thm:ppartdist}, a distribution result on $Q$-partitions on a product of two chains. Note that this is not a new result. Theorem 32 in \cite{GR2015} is a ``reciprocity'' property on rectangular-shaped posets under birational rowmotion which implies Theorem~\ref{thm:ppartdist} via tropicalization. We present a different proof below showing how this result also follows from multifold symmetry.

\begin{theorem}[\protect{\cite[Theorem 32]{GR2015}}]
  \label{thm:ppartdist}
  Let $(x_1,y_1),(x_2,y_2) \in [a] \times [c]$ where $(x_1,y_1)$ and $(x_2,y_2)$ are antipodes and $\sigma \in \mathcal{A}^\ell([a] \times [c])$. Then $\mathrm{Dist}(\sigma,(x_2,y_2),\mathrm{id},v)=\{\ell-m : m \in \mathrm{Dist}(\sigma,(x_1,y_1),\mathrm{id},v) \}$.
\end{theorem}

\begin{proof}
  Our strategy will be to use our bijection $\Phi$ to translate the semistandard Young tableaux distribution result of Theorem~\ref{ssytbcdist} to the desired distribution result on $\sigma \in \mathcal{A}^\ell([a] \times [1] \times [c]) \cong \mathcal{A}^\ell([a] \times [c])$.  We will specialize Corollary \ref{cor:productofthreechains} on $\mathcal{A}^\ell([a] \times [b] \times [c])$ to $b=1$; as part of this we suppress the superfluous middle component in the notation. By Corollary \ref{cor:productofthreechains}, $\Phi\left(\pro(T)\right)=\mathrm{Pro}_{\pi,(-1,1)}\left(\Phi(T)\right)$ where $\pi(i,k)=(i,i-k+c)$. With this choice of $\pi$, $\Phi(\mathcal{L}_{([a] \times [1]) \times [\ell]}(R^{a+c}))$ under $\pro_{\pi, (-1,1)}$ is equivalent to $\mathcal{A}^\ell([a] \times [c])$ under $\pro_{\mathrm{id}, (-1,1)}$; we use this to connect the translation of the semistandard Young tableaux result to the desired distribution result.

  Observe that $\pi$ is a bijection between $\Phi(\mathcal{L}_{([a] \times [1]) \times [\ell]}(R^{a+c}))$ and $\mathcal{A}^\ell([a] \times [c])$. Furthermore, because $\Phi$ is a bijection, for any $\sigma \in  \mathcal{A}^\ell([a] \times [c])$ there exists $T \in \mathcal{L}_{([a] \times [1]) \times [\ell]}(R^{a+c})\cong \mathrm{SSYT}_{a+c}(a \times \ell)$ such that $\Phi(T)=\pi^{-1}(\sigma)$. Observe that $\pi^{-1}(x_1,y_1)=(x_1,x_1-y_1+c)$ and $\pi^{-1}(x_2,y_2)=\pi^{-1}(a+1-x_1,c+1-y_1)=(a+1-x_1,a+y_1-x_1)$. Translating through the bijection $\Phi$, we have
  \[\sigma(x_1,y_1) = \# \{ (x_1,j) \in [a] \times [\ell] :  T(x_1,j) > x_1-y_1+c \}\]
  and
  \[\sigma(x_2,y_2) = \# \{ (a+1-x_1, j) \in [a] \times [\ell] : T(a+1-x_1, j) > a+y_1-x_1  \}.\]

  We now show the desired distribution result for $v=(-1,1)$. Suppose there are $n_{m}$ values of $m \in \mathrm{Dist}(\sigma,(x_1,y_1),\mathrm{id},(-1,1))$. For each of these, there is a $\sigma_{k}\in \{ \mathrm{Pro}^i_{\mathrm{id},(-1,1)}(\sigma) : 0 \le i \le a+b-1\}$, $k \in \{1,\dots,n_m\}$, such that $\sigma_{k}(x_1,y_1)=m$. Translating this through $\Phi$, we have $m=\# \{ (x_1,j) \in [a] \times [\ell] :  T_{k}(x_1,j) > x_1-y_1+c \}$, where $T_{k} \in \{ \mathrm{Pro}^i(T) : 0 \le i \le a+b-1\}$ is the SSYT that corresponds to $\sigma_{k}$ under the bijection $\Phi$. This implies there are $n_m$ values of $m \in \mathrm{BCDist}(T,x_1,x_1-y_1+c)$. By Theorem \ref{ssytbcdist}, there are $n_m$ values of $\ell - m \in \mathrm{BCDist}(T,a+1-x_1,a+y_1-x_1)$. Translating this back through $\Phi$ implies there are $n_m$ values of $\ell - m \in \mathrm{Dist}(\sigma,(x_2,y_2),\mathrm{id},(-1,1))$, completing the proof for $v=(-1,1)$.

  To obtain the result for any $v  \in \{\pm 1\}^2$, we use the recombination technique of Einstein and Propp \cite[Corollary 6.3]{EP2021}, which shows that the distribution of an element in a product of two chains $Q$-partition obtained from a piecewise-linear promotion action is invariant. From this, we obtain the result for any $v  \in \{\pm 1\}^2$.
\end{proof}

\begin{example}
  \label{ex:sstytoppart}
  The first two lines of Figure~\ref{fig:pstricthomomesyexample} show a bijection between $\mathrm{SSYT}_{a+c}(a \times \ell)$ under $\pro$ and $\mathcal{L}_{[a] \times [\ell]}(R^{a+c})$ under $\pro$ for $a=2, c=3, \ell=2$. Reinterpreting $\mathcal{L}_{[a] \times [\ell]}(R^{a+c})$ as  $\mathcal{L}_{([a] \times [1]) \times [\ell]}(R^{a+c})$, lines two and three show our bijection $\Phi$ between $\mathcal{L}_{([a] \times [1]) \times [\ell]}(R^{a+c})$ under $\pro$ and  $\mathcal{A}^\ell([a] \times [1] \times [c]) \cong \mathcal{A}^\ell([a] \times [c])$ under $\pro_{\mathrm{id},(-1,1)}$ for the same parameters. As a result, these three lines show how we translate the distribution result of Theorem~\ref{ssytbcdist} to $Q$-partitions. Lines three and four demonstrate the recombination technique on $\mathcal{A}^\ell([a] \times [c])$ from $\pro_{\mathrm{id},(-1,1)}$ to $\pro_{\mathrm{id},(-1,-1)}$, which we will use in the proof of Theorem~\ref{thm:pstrictdist}. We discuss the connection between lines four and five in Example~\ref{ex:pparttopstrict}.
\end{example}

In Corollary \ref{cor:aborder}, we saw that $\mathcal{L}_{([a] \times [b]) \times [\ell]}(R^{a+b})$ has nice order under $\pro$, and in Corollary~\ref{cor:pstrict_homomesy}, we will obtain a homomesy result on $\mathcal{L}_{([a] \times [b]) \times [\ell]}(R^{a+b})$. These results are somewhat surprising as $[a] \times [b] \times [\ell]$ is an arbitrary product of three chains. The reason we have such nice results is that $\mathcal{L}_{([a] \times [b]) \times [\ell]}(R^{a+b})$ has a rigid structure, with each element having only two possible labels. We make this observation in the following remark.

\begin{remark}
  \label{remark:pstrict_max_aplusb}
  Let $f \in \mathcal{L}_{([a] \times [b]) \times [\ell]}(R^{a+b})$. Then $f(x,y,z)=x+y-1$ or $f(x,y,z)=x+y$ for all $(x,y,z) \in [a] \times [b] \times [\ell]$.
\end{remark}

We now use Theorem~\ref{thm:ppartdist} and Corollary \ref{cor:productofthreechains} to prove a theorem on the distribution of $P$-strict labelings $\mathcal{L}_{([a] \times [b]) \times [\ell]}(R^{a+b})$ under $\pro$.

\begin{theorem}
  \label{thm:pstrictdist}
  Let $(x_1,y_1,z_1),(x_2,y_2,z_2) \in [a] \times [b] \times [\ell]$ where $(x_1,y_1,z_1)$ and $(x_2,y_2,z_2)$ are antipodes and $f \in \mathcal{L}_{([a] \times [b]) \times [\ell]}(R^{a+b})$. Then $\mathrm{Dist}(f,(x_2,y_2,z_2))=\{a+b+1-m : m \in \mathrm{Dist}(f,(x_1,y_1,z_1)) \}$.
\end{theorem}

\begin{proof}
  Our strategy will be to use $\Phi$ to translate the result of Theorem~\ref{thm:ppartdist} to the $P$-strict labelings $\mathcal{L}_{([a] \times [b]) \times [\ell]}(R^{a+b})$. We specialize Corollary~\ref{cor:productofthreechains} to $c=1$: for any $f\in\mathcal{L}_{([a] \times [b]) \times [\ell]}(R^{a+b})$, $\Phi\left(\pro(f)\right)=\mathrm{Pro}_{\pi,(-1,-1,1)}\left(\Phi(f)\right)$ where $\pi((i,j),k)=(i,j,i+j-k)$. Because of the labeling convention from the $\Gamma$ poset construction, the promotion action $\mathrm{Pro}_{\pi,(-1,-1,1)}\left(\Phi(f)\right)$ where $\pi((i,j),k)=(i,j,i+j-k)$ is equivalent to the promotion action $\mathrm{Pro}_{\mathrm{id},(-1,-1)}\left(\sigma\right)$ where $\sigma \in \mathcal{A}^\ell([a] \times [b])$ is the $Q$-partition that corresponds to $\Phi(f)$. As a result, we will eventually apply Theorem~\ref{thm:ppartdist} with this $\sigma$ and $v=(-1,-1)$.

  From a given $f\in\mathcal{L}_{([a] \times [b]) \times [\ell]}(R^{a+b})$, we can determine values of $\sigma$ through the bijection $\Phi$, obtaining
  \[\sigma(x_1,y_1) = \# \{ (x_1,y_1,j) \in [a] \times [b] \times [\ell] :  f(x_1,y_1,j) > x_1+y_1-1 \}\]
  and
  \[\sigma(x_2,y_2) = \# \{ (a+1-x_1,b+1-y_1,j) \in [a] \times [b] \times [\ell] : f(a+1-x_1,b+1-y_1,j) > a+b+1-x_1-y_1 \}.\]

  Recall from Remark \ref{remark:pstrict_max_aplusb} that $f(x,y,z)=x+y-1$ or $f(x,y,z)=x+y$ for all $(x,y,z) \in [a] \times [b] \times [\ell]$. Suppose there are $n$ values of $x_1+y_1-1 \in \mathrm{Dist}(f,(x_1,y_1,z_1))$. For each of these, there is an $f_k \in \{ \mathrm{Pro}^i(f) : 0 \le i \le a+b-1\}$, $k \in \{1,\dots,n\}$, such that $f_k(x_1,y_1,z_1)=x_1+y_1-1$. Translating this through $\Phi$ implies $\sigma_k(x_1,y_1)\le \ell - z_1$, where $\sigma_k \in \{ \mathrm{Pro}^i_{\mathrm{id},(-1,-1)}(\sigma) : 0 \le i \le a+b-1\}$ is the $Q$-partition corresponding to $f_k$ under $\Phi$. By Theorem~\ref{thm:ppartdist}, there must be $n$ values of $\sigma_{k'}\in \{ \mathrm{Pro}^i_{\mathrm{id},(-1,-1)}(\sigma) : 0 \le i \le a+b-1\}$ such that $\sigma_{k'}(a+1-x_1,b+1-y_1) \ge z_1$. As a result, there must be at least $n$ values of $f_{k'}\in \{ \mathrm{Pro}^i(f) : 0 \le i \le a+b-1\}$ such that $f_{k'}(a+1-x_1,b+1-y_1,\ell+1-z_1)=a+b+2-x_1-y_1$. Ergo, there must be at least $n$ values of $a+b+1-(x_1+y_1-1) \in \mathrm{Dist}(f,(x_2,y_2,z_2))$.

  We now consider the $a+b-n$ values of $x_1+y_1 \in \mathrm{Dist}(f,(x_1,y_1,z_1))$. For each of these, there is an $f_K \in \{ \mathrm{Pro}^i(f) : 0 \le i \le a+b-1\}$, $K \in \{1,\dots,a+b-n\}$, such that $f_K(x_1,y_1,z_1)=x_1+y_1$. Translating this through $\Phi$ implies $\sigma_K(x_1,y_1) > \ell - z_1$ where $\sigma_K \in \{ \mathrm{Pro}^i_{\mathrm{id},(-1,-1)}(\sigma) : 0 \le i \le a+b-1\}$ is the $P$-partition corresponding to $f_K$ under $\Phi$. Once again, we apply Theorem~\ref{thm:ppartdist} to obtain $a+b-n$ values of $\sigma_{K'}\in \{ \mathrm{Pro}^i_{\mathrm{id},(-1,-1)}(\sigma) : 0 \le i \le a+b-1\}$ such that $\sigma_{K'}(a+1-x_1,b+1-y_1) < z_1$. As a result, there must be at least $a+b-n$ values of $f_{K'}\in \{ \mathrm{Pro}^i(f) : 0 \le i \le a+b-1\}$ such that $f_{K'}(a+1-x_1,b+1-y_1,\ell+1-z_1)=a+b+1-x_1-y_1$. From this, there must be at least $a+b-n$ values of $a+b+1-(x_1+y_1)\in \mathrm{Dist}(f,(x_2,y_2,z_2))$.

  Because there are at least $n$ values of $a+b+1-(x_1+y_1-1) \in \mathrm{Dist}(f,(x_2,y_2,z_2))$, at least $a+b-n$ values of $a+b+1-(x_1+y_1)\in \mathrm{Dist}(f,(x_2,y_2,z_2))$, and only $a+b$ total values in $\mathrm{Dist}(f,(x_2,y_2,z_2))$, this gives us exactly $n$ values of $a+b+1-(x_1+y_1-1) \in \mathrm{Dist}(f,(x_2,y_2,z_2))$ and exactly $a+b-n$ values of $a+b+1-(x_1+y_1)\in \mathrm{Dist}(f,(x_2,y_2,z_2))$, which completes the proof.
\end{proof}

\begin{corollary}
  \label{cor:pstrict_homomesy}
  Let $S=\{x,y\} \subseteq [a] \times [b] \times [\ell]$ where $x$ and $y$ are antipodes. Then \[(\mathcal{L}_{([a] \times [b]) \times [\ell]}(R^{a+b}), \mathrm{Pro}, \chi_S)\] is $(a+b+1)$-mesic.
\end{corollary}

\begin{proof}
  From Theorem~\ref{thm:pstrictdist}, we know $\mathrm{Dist}(f,x)=\{a+b+1-i : i \in \mathrm{Dist}(f,y) \}$. As a result, if we average the statistic $\chi_S$ over an orbit $\mathcal{O}$ of $\mathcal{L}_{([a] \times [b]) \times [\ell]}(R^{a+b})$ under $\pro$, we obtain
  \[ \displaystyle\frac{1}{\# \mathcal{O}}\left(\sum_{i \in \mathrm{Dist}(f,x)} i+\sum_{i \in \mathrm{Dist}(f,y)}i\right)=\frac{1}{\# \mathcal{O}}\sum_{i=1}^{\# \mathcal{O}}[a+b+1-i+i]=\frac{\# \mathcal{O}(a+b+1)}{\# \mathcal{O}}=a+b+1. \]
\end{proof}

\begin{example}
  \label{ex:pparttopstrict}
  Lines four and five of Figure~\ref{fig:pstricthomomesyexample} show our bijection $\Phi$ between $\mathcal{A}^\ell([a] \times [b])$ under $\pro_{\mathrm{id},(-1,-1)}$ and $\mathcal{L}_{([a] \times [b]) \times [\ell]}(R^{a+b})$ under $\pro$ where $a=2, b=3, \ell=2$. The entire figure shows our strategy of the proofs of Theorems~\ref{thm:ppartdist} and \ref{thm:pstrictdist}, translating the distribution result of semistandard Young tableaux under $\pro$ to $\mathcal{L}_{([a] \times [b]) \times [\ell]}(R^{a+b})$ under $\pro$.
\end{example}

\subsection{Conjectures}
\label{subsec:conj}
Based on computational evidence, we conjecture the following antipodal symmetry and translate this conjecture using our main bijection and multifold symmetry.
\begin{conjecture} Let $S=\{x,y\} \subseteq [a] \times [b] \times [\ell]$ where $x$ and $y$ are antipodes. Then \[\left(\mathcal{A}^2([a] \times [2] \times [2]),\row,\chi_S\right)\] is $2$-mesic.
\end{conjecture}

We have verified this conjecture holds for $a\leq 6$.

\smallskip
By Corollary~\ref{cor:productofthreechains} and by translating $\chi$ through $\Phi$, this is equivalent to the following conjecture.

\begin{definition}
  Let $f \in \mathcal{L}_{([a] \times [2]) \times [2]}(R^{a+3})$. Define the statistic $\xi(f,(x_1,x_2),b)$ to be
  \[\xi(f,(x_1,x_2),b) = \displaystyle\sum_{k=1}^b k \cdot \# \{ (x_1,x_2,j) \in [a] \times [2] \times [2] :  f(x_1,x_2,j) = x_1+x_2+k-1 \}.\]
\end{definition}

\begin{conjecture}
  Let $a \in \mathbb{N}$. Then
  \[(\mathcal{L}_{([a] \times [2]) \times [2]}(R^{a+3}), \pro, \xi(f,(x_1,x_2),2)+\xi(f,(a+1-x_1,3-x_2),2)) \] is $4$-mesic.
\end{conjecture}

By Theorem~\ref{thm:3dmultifold} and by translating $\chi$ through $\Phi$, we have the following equivalent conjecture.

\begin{conjecture}
  The triple
  \[(\mathcal{L}_{([2] \times [2]) \times [2]}(R^{a+3}), \pro, \xi(f,(x_1,x_2),a)+\xi(f,(3-x_1,3-x_2),a)) \]
  is $4$-mesic.
\end{conjecture}

\section{Beyond the product of chains}
\label{sec:beyond}
The previous section studied $P$-strict promotion where $P$ is a product of chains and the restriction function is induced by a global bound. In this section, we apply our main theorem to examples of interest where $P$ is not a product of chains or the restriction function is not induced by a global bound. In Subsection~\ref{subsec:minuscule}, we study the case of $P$ a minuscule poset. In Subsection~\ref{subsec:flagged}, we let $P=[a]\times[b]$ but impose flags on our restriction function. Finally, Subsection~\ref{sec:vee} discusses the case where $P$ is the three-element poset $V$.

\subsection{Minuscule posets}
\label{subsec:minuscule}
Minuscule posets are interesting families of posets arising from Lie theory. The product of two chains $[a]\times[b]$ is a \emph{Type $A$ minuscule poset}, and thus, it is natural to ask which results of the previous section extend to other types. The answer is: not many. A special feature of the Type $A$ minuscule poset is that it is constructed as a Cartesian product of chains, and many of the results in Section~\ref{sec:ab} rely on the multifold symmetry of Corollary~\ref{thm:3dmultifold}. Since the other minuscule posets are not products of chains, we do not have an analogous result to use. But, our main bijection is useful for translating results on $Q$-partitions where $Q$ is the Cartesian product of a minuscule poset and a chain to obtain new results on $P$-strict labelings of minuscule posets. The new such result of this section is Corollary~\ref{cor:a1}.

Since the focus of this section is on translating known results about minuscule posets across our bijection, we choose to give explicit descriptions of the minuscule posets rather than the full Lie-theoretic definitions. Rather, we refer the reader to the cited papers for further explanation of the algebraic meaning.

The minuscule posets are the following three infinite families followed by two exceptional posets. We follow the notation and definitions of minuscule posets found in \cite{PechenikMinuscule}. For each, we also give the \emph{Coxeter number} ${h}$ of the associated Lie algebra, which will be used in statements of results.
Let $J^k$ denote the order ideal functor applied $k$ times. For example, $J^2(P)=J(J(P))$ is the poset of order ideals of the poset of order ideals of $P$.
\begin{enumerate}
  \item Rectangles: $[k]\times[m]$, $h=k+m$,
  \item Shifted staircases: $S_k:=\{(x,y) \ | \ x\leq y\in [k]\}$, $h=2k$
  \item Propellers: $J^k([2]\times[2])$, $h=2(k+2)$,
  \item Cayley--Moufang: $J^2([3]\times[2])$, $h=12$
  \item Freudenthal: $J^3([3]\times[2])$, $h=18$.
\end{enumerate}

First, since all minuscule posets are graded, we can use Corollary~\ref{cor:PrankedGlobalqRow} to construct an equivariant bijection.
\begin{corollary}
  \label{cor:minuscule}
  Let $P$ be a minuscule poset of rank $n$. Then $\mathcal{A}^{\ell}(P\times[a])$ under $\row$ is in equivariant bijection with $\mathcal{L}_{P\times [\ell]}(R^{a+n+1})$ under $\pro$.
\end{corollary}
\begin{proof}
  This follows directly from Corollary~\ref{cor:PrankedGlobalqRow} and the fact that any minuscule poset is graded.
\end{proof}

We state some known results on $\mathcal{A}^{\ell}(P\times[a])$ and/or $\mathcal{L}_{P\times [\ell]}(R^{a+n+1})$ and their translations via Corollary~\ref{cor:minuscule}.

In the case $\ell=1$, the bijection of Corollary~\ref{cor:minuscule} was already known and used by H.~Mandel and Pechenik  to obtain results on order and cyclic sieving in minuscule posets. We state these below, together with previous results of D.~Rush and X.~Shi which used the order ideal perspective. (For more on the $\ell=1$ case for general $P$, see also Remark~\ref{remark:ell1}.)

\emph{Cyclic sieving} is a dynamical phenomenon in which evaluation of a polynomial at certain roots of unity completely describes the orbit structure of an action. Since we state no new cyclic sieving results here, we refer the reader to the foundational paper~\cite{ReStWh2004} and the papers cited in the theorem below for the precise definition.

\begin{theorem}[\protect{\cite[$a=1,2$ on order ideals]{RS2013}}, \protect{\cite[other cases]{MandelPechenik}}]
  Let $P$ be a minuscule poset associated with a minuscule weight $\lambda$ of $\mathfrak{g}$, and let $n$ denote the rank of $P$. Each of
  \begin{itemize}
    \item the set of order ideals $J(P\times [a])$  under $\row$ and
    \item the set of increasing labelings $\mathcal{L}_{P\times [1]}(R^{a+n+1})$ under $\pro$
  \end{itemize}
  have order dividing $h$ and exhibit the cyclic sieving phenomenon with respect to the cardinality generating function of $J(P\times [a])$ for the values of $a$ given below:
  \begin{enumerate}
    \item Rectangles: $a=1,2$,
    \item Shifted staircases: $a=1,2$,
    \item Propellers: all $a$,
    \item Cayley--Moufang: all $a$,
    \item Freudenthal: $a\leq 4$.
  \end{enumerate}
\end{theorem}

The next theorem and translation allow for arbitrary $\ell$, but specify $a=1$. The theorem has appeared in the literature in several places; for more on these references, see the discussion in \cite[p.\ 3-5]{Okada2021}.
\begin{theorem}[\cite{GarverPatriasThomas,GR2015,GR2016,Okada2021}]
  \label{thm:PpartMinuscule}
  Let $P$ be a minuscule poset associated with a minuscule weight $\lambda$ of $\mathfrak{g}$, and let $n$ denote the rank of $P$. Then $\mathcal{A}^{\ell}(P)$ has order $h$ under $\row$, where $h$ is the Coxeter number of $\mathfrak{g}$.
\end{theorem}

The following corollary is a translation of the above theorem, using our main bijection. The Type A case was previously stated as Corollary~\ref{cor:aborder}.
\begin{corollary}
  \label{cor:a1}
  Let $P$ be a minuscule poset associated with a minuscule weight $\lambda$ of $\mathfrak{g}$, and let $n$ denote the rank of $P$. Then
  the set of $P$-strict labelings $\mathcal{L}_{P\times [\ell]}(R^{n+2})$ has order $h$ under $\mathrm{Pro}$, where $h$ is the Coxeter number of $\mathfrak{g}$.
\end{corollary}
\begin{proof}
  By Corollary~\ref{cor:minuscule}, $\mathcal{A}^{\ell}(P\times[1])$ is in bijection with $\mathcal{L}_{P\times [\ell]}(R^{n+2})$. Then the result follows by Theorem~\ref{thm:PpartMinuscule}.
\end{proof}

In addition, \cite{Okada2021} includes discussion of several homomesy results, which one could translate, but this would require more background on minuscule posets than what we have stated in this section.

One might hope for some nice results on the order of rowmotion on $\mathcal{A}^{\ell}(P\times[a])$ when both $a$ and $\ell$ are greater than one (and $P$ is a minuscule poset). But our SageMath computations show that for the propeller $P=J^2([2]\times[2])$, $\mathcal{A}^{2}(P\times[2])$ has rowmotion orbits of sizes $6$, $9$, $17$, and $44$. Since this is one of the simplest minuscule posets, a nice general result involving order or cyclic sieving is unlikely.

\subsection{Flagged product of two chains}
\label{subsec:flagged}
A \emph{flagged tableau} of $n$ rows with flag $\beta = (\beta_1,\beta_2,\ldots,\beta_n)$ is a semistandard tableau whose entries in row $i$ do not exceed $\beta_i$.  Given that semistandard tableaux are $P$-strict labelings with $P = [n]$, it is natural to consider a flagged $P$-strict labeling with $P \neq [n]$, where the flags on the fibers (the ``rows" of our $P$-strict labeling) are given by the function $\beta: P \to \mathbb{Z}^+$.  In this section, we explore a case where a particular flagging of $[a] \times [b]$ results in a familiar structure for the gamma poset.

The poset we consider is the following.
\begin{definition}
  Let $\widetriangle_n$ denote the triangle-shaped subposet of $[n] \times [n]$ given by $\{(i,j) \mid 1 \leq i \leq n, n-i+1 \leq j \leq n\}$.
\end{definition}

\begin{theorem} \label{thm:atimesbtypea}
  For $(i,j) \in [a] \times [b]$, let $\beta(i,j) = b + 2i - 1$.  Then $\mathcal{L}_{([a] \times [b]) \times [\ell]}(R^\beta)$ under $\pro$ is in equivariant bijection with $\mathcal{L}_{\subwidetriangle_{a} \times [\ell]}(R^{a+b})$ under $\pro$.
\end{theorem}

\begin{figure}
  \begin{center}
    \includegraphics{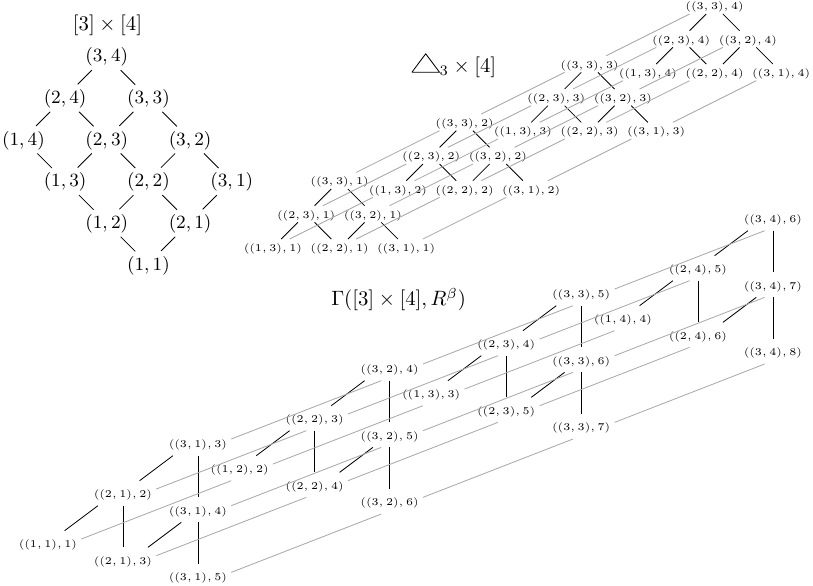}
  \end{center}
  \caption{An example of the posets in Lemma \ref{lem:atimesbtypea} with $a = 3$ and $b = 4$. The restriction function $R^\beta$ on $[3] \times [4]$ is defined by the flags $\beta$, where $\beta(1,j) = 5$, $\beta(2,j) = 7$, and $\beta(3,j) = 9$.  Note that the flags do not necessarily give the greatest elements of the consistent restriction function associated to $R^\beta$ because of the strictly increasing requirement on layers in $\mathcal{L}_{([a] \times [b]) \times [\ell]}(R^\beta)$.} \label{fig:atimesbtypea}
\end{figure}

To prove Theorem \ref{thm:atimesbtypea}, we use the following lemma, which proves the poset isomorphism shown in Figure \ref{fig:atimesbtypea}.

\begin{lemma} \label{lem:atimesbtypea}
  Let $\beta$ be defined as above. Then $\Gamma([a] \times [b], R^\beta)$ is isomorphic as a poset to $\widetriangle_{a} \times [b]$.
\end{lemma}

\begin{proof}
  The consistent restriction function $R^B_A$ associated to $R^{\beta}$ is  \[R^B_A(i,j) = \{i+j-1, i+j,\ldots,  b+2i-1-(b-j)\} = \{i+j-1,\ldots,2i+j-1\},\] so $\Gamma([a] \times [b], R^\beta)$ consists of the elements \[\{((i,j),k) \mid 1 \leq i \leq a, 1 \leq j \leq b, i+j-1 \leq k \leq 2i+j-2\}.\]  By definition of $\Gamma$, $((i_1,j_1), k_1) \lessdot ((i_2,j_2),k_2)$ in $\Gamma([a] \times [b], R^\beta)$ if and only if either \begin{enumerate}
    \item[(1)] $(i_1,j_1) = (i_2,j_2)$ and $k_1-1= k_2$, or
    \item[(2)] $k_1+1= k_2$ and $(i_1,j_1) \lessdot_{[a]\times [b]} (i_2,j_2)$.
  \end{enumerate}
  Consider the map from $\Gamma([a] \times [b], R^\beta)$ to $\widetriangle_{a} \times [b]$ defined by \[((i,j),k) \mapsto ((i,i+j-k+a-1),j) := ((i,u),j).\]  Then, using the bounds on $k$, we have
  \[i+j-(2i+j-2)+a-1 \leq u \leq i+j-(i+j-1)+a-1,\]
  implying $a-i+1\leq u \leq a$. Together with the bounds $1 \leq i \leq a$ and $1 \leq j \leq b$, we conclude our map is a bijection between the elements of $\Gamma([a] \times [b], R^\beta)$ and $\widetriangle_{a} \times [b]$.  Moreover, 
  $((i_1,u_1),j_1) \lessdot ((i_2,u_2),j_2)$ in ${\widetriangle_a \times [b]}$
  if and only if either $i_1+1 = i_2$, $u_1+1 = u_2$, or $j_1+1 = j_2$, and the other coordinates are equal.  Each of these conditions correspond exactly with the covering relations in $\Gamma([a] \times [b], R^\beta)$ as follows:  if $i_1+1 = i_2$ or $j_1+1 = j_2$, then $(i_1,j_1) \lessdot_{[a]\times [b]} (i_2,j_2)$ and, since $u_1 = i_1+j_1-k_1+a-1 = u_2 = i_2+j_2-k_2+a-1$, it must be the case that $k_1 + 1 = k_2$.  If $u_1+1 = u_2$ and $i_1 = i_2,$ $j_1 = j_2$, then $-k_1-1 + 1 = -k_2-1$, so $k_1 - 1 = k_2$.  Thus, the above map bijects the elements of $\Gamma([a] \times [b], R^\beta)$ to those of $\widetriangle_{a} \times [b]$ and preserves covering relations, so $\Gamma([a] \times [b], R^\beta)$ is isomorphic as a poset to $\widetriangle_{a} \times [b]$.
\end{proof}

\begin{proof}[Proof of Theorem \ref{thm:atimesbtypea}]
  Since $\mathcal{L}_{([a] \times [b]) \times [\ell]}(R^\beta)$ under $\pro$ is in equivariant bijection with $\mathcal{A}^\ell(\Gamma([a] \times [b], R^\beta))$ under $\row$ by Corollary \ref{cor:abrow},  we have that it is also in equivariant bijection with $\mathcal{A}^\ell(\widetriangle_{a} \times [b])$ under $\row$ by Lemma \ref{lem:atimesbtypea}.  Since $\widetriangle_{a}$ has rank $a-1$, $\mathcal{A}^\ell(\widetriangle_{a} \times [b])$ under $\row$ is in equivariant bijection with $\mathcal{L}_{\subwidetriangle_{a} \times [\ell]}(R^{a+b})$ under $\pro$ by Corollary~\ref{thm:PrankedGlobalq}, and we have the desired result.
\end{proof}

\begin{remark}
  Lemma \ref{lem:atimesbtypea} appears as  \cite[Lemma 4.27]{BSV2021} for the case $b=1$.
\end{remark}

\begin{remark}
  There is an explicit bijection between $\mathcal{L}_{([a] \times [b]) \times [\ell]}(R^\beta)$ and $\mathcal{L}_{\widetriangle_{a} \times [\ell]}(R^{a+b})$ using $\Phi$ from Definition \ref{def:mainbijection} with the appropriate restriction functions. However, this is not the same as the equivariant bijection from the above theorem, since the toggle orders associated to their respective gamma posets are different (but both conjugate to rowmotion).
\end{remark}

\begin{remark}
  As a result of D.~Grinberg and T.~Roby \cite[Corollary 66]{GR2015}, we have that the order of rowmotion on $\mathcal{A}^\ell(\widetriangle_{a} \times [1])$ is given by $2(a+1)$, and it is conjectured (see Conjecture~\ref{conj:Vtimesm}) that $\mathcal{A}^\ell(\widetriangle_{2} \times [b])$ also has predictable order under rowmotion.  Additionally, based on some preliminary calculations in SageMath, there may be a resonance result (see Definition~\ref{def:resonance}) on the order ideals $\mathcal{A}^1(\widetriangle_{a} \times [b])$. Unfortunately, any change to the parameters beyond those of the sets listed above results in a breakdown of the orbit structure.  As early as $\mathcal{A}^2(\widetriangle_{3} \times [2])$, we find an orbit of length $94$. It could be said, however, that \emph{many} of the orbits are of length $10$, a notion which may be worth future consideration.
\end{remark}

\subsection{The V poset}
\label{sec:vee}
Corollaries~\ref{thm:PrankedGlobalq} and \ref{cor:PrankedGlobalqRow} give us reason to pursue any cases where $P$ is graded and $\row$ on $\mathcal{A}^\ell(P \times [m])$ has noteworthy dynamical properties, as the translation through our bijection produces an example in which $\pro$ on $\mathcal{L}_{P \times [\ell]}(R^q)$ is also of dynamical interest.  One such poset that has caught recent attention is $P=V \times [m]$, where $V$ is defined below.

\begin{definition}
  Let $V$ be the three-element poset $\{a,b,c\}$ with $a \lessdot b$ and $a \lessdot c$.
\end{definition}

As a result of M.~Plante, we know the order ideals of $V \times [m]$ have nice order under rowmotion.

\begin{theorem}[\cite{Plante2022}]
  $\row$ on $J(V \times [m])$ has order dividing $2(m+2)$.
\end{theorem}

Moreover, it is believed to be the case that piecewise-linear rowmotion on $(V\times[m])$-partitions shares the same order.

\begin{conjecture}[\cite{hopkins2021order}] \label{conj:Vtimesm}
  $\row$ on $\mathcal{A}^\ell(V \times [m])$ has order dividing $2(m+2)$.
\end{conjecture}

This conjecture has been experimentally verified for small $m$ and $\ell$, and we have used our SageMath code to check many particular examples with larger $\ell$ and $m$.

We apply Corollary~\ref{thm:PrankedGlobalq}  to Conjecture \ref{conj:Vtimesm} (noting that the poset $V$ is of rank $n = 1$) to obtain the following translated conjecture on $V$-strict labelings:

\begin{conjecture} \label{conj:Vtimesell}
  $\pro$ on $\mathcal{L}_{V \times [\ell]}(R^q)$ is of order dividing $2q$.
\end{conjecture}

\begin{remark}
  Hopkins and M.~Rubey show in \cite{HopkinsRubey2021} that the order of promotion on linear extensions of $V \times [\ell]$ is $6\ell$. Just as in the linear extension case, it seems to be true (experimentally) that applying promotion $q$ times to some $f \in \mathcal{L}_{V \times [\ell]}(R^{q})$ results in a reflection of the labels of $f$ across the vertical axis of symmetry. Therefore, it may be possible to devise a method similar to that in \cite{HopkinsRubey2021} to accommodate $V$-strict labelings, proving Conjecture~\ref{conj:Vtimesell}, and, in turn, Conjecture~\ref{conj:Vtimesm}.
\end{remark}

Finally, Conjecture \ref{conj:Vtimesm} can be translated once more to an earlier case of interest, by applying Theorem \ref{thm:atimesbtypea} with $a=2$, noting $\widetriangle_2$ is the dual poset of $V$.

\begin{conjecture} For $(i,j) \in [2] \times [b]$, let $\beta(i,j) = b + 2i-1$.  Then $\pro$ on $\mathcal{L}_{([2] \times [b]) \times [\ell]}(R^\beta)$ is of order dividing $2(b+2)$.
\end{conjecture}

\section{Resonance on P-strict labelings}
\label{sec:resonance}
In dynamical algebraic combinatorics, we are often interested in when an action, such as promotion or rowmotion, has a small, predictable order. Actions which do not have such an order may still exhibit some nice dynamical behavior, such as orbit sizes that may be multiples of a predictable number (or divisors of it). In \cite{DPS2017}, \emph{resonance} was defined to explain this numerical phenomenon when the action in question projects to an action with a small, predictable order.
\begin{definition}[\cite{DPS2017}] \label{def:resonance}
  Suppose $G = \langle g\rangle$ is a cyclic group acting on a set $X$, $C_{\omega} = \langle c\rangle$ a cyclic group of
  order $\omega$ acting nontrivially on a set $Y$, and $\varphi:X \rightarrow Y$ a surjection. We say the triple $(X, G, \varphi)$
  exhibits \textbf{resonance} with \textbf{frequency} $\omega$ if, for all $x\in X$, $c \cdot \varphi(x) = \varphi(g \cdot x)$.
\end{definition}

A prototypical example of resonance given in \cite{DPS2017} was increasing tableaux under K-promotion where the projection map $\varphi$ was the binary content of the tableau. Here we give an analogue of that theorem in the more general setting of $P$-strict labelings.

Though the previous sections of this paper deal only with the case $\mathcal{L}_{P\times[\ell]}(R^q)$, which in the notation of \cite{BSV2021} is $\mathcal{L}_{P\times[\ell]}(u,v,R^q)$ with $u=v=0$, the theorems of this section apply to the case with general $u$ and $v$. As the values of $u$ and $v$ play no role in the proof, we state the first few results of this section at this greater level of generality. See~\cite[Definition 1.7]{BSV2021} for the definition of $\mathcal{L}_{P\times[\ell]}(u,v,R^q)$.

We first extend the definition of binary content from increasing tableaux to $P$-strict labelings.
\begin{definition}
  Define the \textbf{binary content} of a $P$-strict labeling $f \in \mathcal{L}_{P\times[\ell]}(u,v,R^q)$
  to be the sequence
  $\operatorname{Con}(f) = (a_1, a_2,\ldots , a_q)$, where $a_i = 1$ if $f(p,i) = i$ for some $(p,i) \in P \times [\ell]$ and $0$ otherwise.
\end{definition}

We now give Lemma~\ref{lem:cyclic_shift} showing that promotion cyclically shifts the binary content of a $P$-strict labeling. This is an analogue of \cite[Lemma 2.1]{DPS2017}. Note the proof of Lemma~\ref{lem:cyclic_shift} is not directly analogous to the proof of \cite[Lemma 2.1]{DPS2017}, as we use the Bender--Knuth definition of $P$-strict promotion (Definition~\ref{def:BenderKnuth}) rather than an analogue of \textit{jeu de taquin}. In \cite[Definition 3.1]{BSV2021}, we gave a definition of $P$-strict promotion via an analogue of jeu de taquin when the restriction function is $R^q$ and showed these definitions are indeed equivalent \cite[Theorem 3.10]{BSV2021}. So we could have used this approach instead, at the cost of stating this alternate definition for $P$-strict promotion.
\begin{lemma}
  \label{lem:cyclic_shift}
  Let $f \in \mathcal{L}_{P\times[\ell]}(u,v,R^q)$. If $\operatorname{Con}(f) = (a_1, a_2,\ldots , a_q)$, then $\operatorname{Con}(\pro(f))$ is the cyclic shift
  $(a_2,\ldots , a_q, a_1)$.
\end{lemma}
\begin{proof}
  Suppose $f\in \mathcal{L}_{P\times[\ell]}(u,v,R^q)$, $\operatorname{Con}(f) = (a_1, a_2,\ldots , a_q)$, and $\operatorname{Con}(\pro(f)) = (b_1, b_2,\ldots , b_q)$. We wish to show $b_i=a_{i+1}$ for all $1\leq i< q$ and $b_q=a_1$.

  Suppose $a_1=0$, meaning $1$ is not used as a label in $f$. Then all $2$ labels in $f$ will be lowerable, leaving no $2$ labels after the application of $\rho_1$. Likewise, there are no $i+1$ labels in $\rho_i\rho_{i-1}\cdots\rho_1(f)$, and there are no $q$ labels in $\pro(f)=\rho_{q-1}\cdots\rho_2\rho_1$. Thus $q$ is not used as a label in $\pro(f)$, therefore $b_q=0=a_1$.

  Suppose $a_1=1$, meaning $1$ is used as a label in $f$. Any $1$ label is either fixed, meaning there is a $2$ above it in its layer, or free. If it is free, then either there are no $2$'s in its fiber and the $1$ changes to a $2$, or there are $2$'s in the fiber, meaning there will still be at least one $2$ in the fiber after the application of $\rho_1$. The same reasoning holds for $i$, $2\leq i\leq q-1$. So there will be at least one $q$ in $\pro(f)$, therefore $b_q=1=a_1$.

  Likewise, suppose $i>1$ and $a_i=0$, meaning $i$ is not used as a label in $f$. Then all $i-1$ labels in $f$ will be raisable, leaving no $i-1$ labels in $\rho_i\rho_{i-1}\cdots\rho_1(f)$, and therefore in $\pro(f)$. So $b_{i-1}=0=a_i$.

  Finally, suppose $i>1$ and $a_i=1$, so there is at least one $i$ used as a label in $f$. Thus, $i$ is fixed in $\rho_{i-2}\cdots\rho_1(f)$ with respect to $\rho_{i-1}$, meaning there is an $i-1$ below it in its layer. In that case, there is an $i-1$ in $\pro(f)$. Alternatively, the $i$ is free, meaning it will either change into an $i-1$ when $\rho_i$ is applied, or there will be $a$ labels of $i-1$ and $b$ labels of $i$ that change to $b$ labels of $i-1$ and $a$ labels of $i$. In either case, there is an $i-1$ in $\pro(f)$. As a result, $b_{i-1}=1=a_i$.
\end{proof}

\begin{theorem} \label{thm:ConResonance}
  $(\mathcal{L}_{P\times[\ell]}(u,v,R^q),\pro,\operatorname{Con})$ exhibits resonance with frequency $q$.
\end{theorem}
\begin{proof}
  This follows directly from Lemma~\ref{lem:cyclic_shift}.
\end{proof}

We now translate this resonance theorem to the realm of $Q$-partitions via the bijections of Theorem~\ref{thm:moregeneralpro} and Corollary~\ref{cor:abrow}. Though we could state this result for arbitrary $u,v$ through the general bijection of \cite{BSV2021}, it would use too much excess notation. Thus, we restrict to the case $u=v=0$. First, we define the analogue of $\operatorname{Con}$ in this realm.

\begin{definition}
  For $\sigma \in \mathcal{A}^{\ell}({\Gamma}(P,{R^q}))$, let $\operatorname{Diff}(\sigma)$ be the sequence $(a_1,a_2,\ldots,a_q)$ where \[a_k = \begin{cases} 0 & \text{if }\sigma(p,k-1) = \sigma(p,k) \text{ for all } p \in P \\
      1 & \text{otherwise,}
    \end{cases}\] where we consider $\sigma(p,j) = \ell$ if $j < \min R^q(p)$ and $\sigma(p,j) = 0$ if $j > \max R^q(p)$.
\end{definition}

\begin{corollary} \label{cor:generalDiff}
  $(\mathcal{A}^{\ell}({\Gamma}(P,{R^q})),\row,\operatorname{Diff})$ exhibits resonance with frequency $q$.
\end{corollary}

\begin{proof}
  This follows from Corollary~\ref{cor:abrow}, Theorem~\ref{thm:ConResonance}, and the fact that, for any $\sigma \in \mathcal{A}^{\ell}({\Gamma}(P,{R^q}))$, the difference $\sigma(p,k-1) - \sigma(p,k)$ gives the number of $k$ labels in the fiber $F_p$ for the corresponding $P$-strict labeling $f$.
\end{proof}

In the case where $P$ is graded, that is, where the poset ${\Gamma}(P,{R^q})$ is isomorphic to $P \times [q-n-1]$, we can restate Corollary \ref{cor:generalDiff} in terms of the elements of $P \times [q-n-1]$ through a reinterpretation of $\operatorname{Diff}$.

\begin{definition}
  Let $P$ be graded with rank $n$ and let $H_k$ be the set of all elements $(p,i)$ in $P \times \mathbb{Z}$ with $q - n - i + \operatorname{rank}(p) = k$.  For $\sigma \in \mathcal{A}^\ell(P \times [q-n-1])$, let $\operatorname{Diff}(\sigma)$ be the sequence $(a_1,a_2,\ldots,a_q)$ where \[a_k = \begin{cases} 0 & \text{if }\sigma(p,i) = \sigma(p,i+1) \text{ for all } (p,i) \in H_{k} \\
      1 & \text{otherwise,}
    \end{cases}\] where we consider $\sigma(p,i) = \ell$ if $i > q-n-1$ and $\sigma(p,i) = 0$ if $i < 1$.
\end{definition}

\begin{corollary}
  $(\mathcal{A}^{\ell}(P \times [q-n-1]),\row,\operatorname{Diff})$ exhibits resonance with frequency $q$.
\end{corollary}
\begin{proof}
  We have that ${\Gamma}(P,{R^q}))$ is isomorphic to $P \times [q-n-1]$ as a poset by Lemma \ref{lem:GradedGlobalq}.  The conditions for $a_k = 0$ and $a_k = 1$ are identical to those in Corollary~\ref{cor:generalDiff} under the bijection in Remark~\ref{rem:proof_bij}, noting that $H_{k'}$ in $P \times [q-n-1]$ corresponds to the set $\{(p,k) \in {\Gamma}(P,{R^q}) \mid k = k'\}$.
\end{proof}

\begin{remark}
  In the case when $P$ is graded, we can think of $\operatorname{Diff}(\sigma)$ as indicating when the elements of consecutive hyperplanes in a lattice projection of $P \times [q-n-1]$ share the same labels.  For example, in the top right of Figure~\ref{fig:PrankedGlobalq}, the two 4 labels at the top of $P \times \{3\}$, the center 1 label in $P \times \{2\}$, and the two 0 labels in $P \times \{1\}$ label the elements of $H_3$, and these labels are the same as those of $H_4$ (where the elements ``below" the poset are considered to be labeled by 0).  In this case, we have $\operatorname{Diff}(\sigma) = (1,1,1,0,1,1)$.
\end{remark}

\section*{Acknowledgments}
The authors thank the referees for their detailed reading and Matthew Plante for helpful conversations regarding Section~\ref{sec:vee}. They also thank the developers of \verb|SageMath|~\cite{sage} software, which was useful in this research, and the \verb|CoCalc|~\cite{SMC} collaboration platform. They thank O.\ Cheong for developing Ipe~\cite{ipe}, which we used to create the figures in this paper. JS was supported by a grant from the Simons Foundation/SFARI (527204, JS) and NSF grant DMS-2247089.

\bibliographystyle{abbrv}
\bibliography{master}
\end{document}